\RequirePackage{xr-hyper}
\documentclass[onefignum,onetabnum]{siamart251216}

\usepackage{lipsum}
\usepackage{amsfonts}
\usepackage{graphicx}
\usepackage{epstopdf}
\usepackage[caption=false]{subfig} 
\usepackage{amssymb}
\usepackage{amsfonts}
\usepackage{arydshln}
\usepackage{mathtools}
\usepackage{algorithm}
\usepackage{algpseudocode}
\usepackage{bm}
\usepackage{comment}
\usepackage[framemethod=tikz]{mdframed}
\ifpdf
\usepackage{xr-hyper}

  \DeclareGraphicsExtensions{.eps,.pdf,.png,.jpg}
\else
  \DeclareGraphicsExtensions{.eps}
\fi

\newsiamremark{remark}{Remark}
\newsiamremark{hypothesis}{Hypothesis}
\crefname{hypothesis}{Hypothesis}{Hypotheses}
\newsiamthm{claim}{Claim}
\newsiamremark{fact}{Fact}
\crefname{fact}{Fact}{Facts}
\newtheorem{assumption}[theorem]{Assumption}
\crefname{assumption}{Assumption}{Assumptions}

\headers{Multi-trajectory parallel-in-iteration optimization}{G. H. M. Ara\'ujo, O. A. Krzysik, and H. de Sterck}

\title{Multi-trajectory parallel-in-iteration optimization}

\author{Guilherme H. M. Ara\'ujo\thanks{Department of Applied Mathematics, University of Waterloo, Waterloo, ON 
  (\email{ghmaciei@uwaterloo.ca}, \email{hdesterck@uwaterloo.ca}).}
\and Oliver A. Krzysik\thanks{Los Alamos National Laboratory, Los Alamos, NM 
  (\email{okrzysik@lanl.gov}).}
\and Hans De Sterck\footnotemark[1]}

\usepackage{amsopn}

\begin{document}

\maketitle
% \begin{abstract}
%     Previously we developed a  ``parallel-in-iteration'' optimization framework that parallelizes across iterations or trajectories of slowly converging sequential optimization methods (e.g., gradient descent) using multigrid reduction-in-time (MGRIT). Here, we develop an extension of this framework.
%     %
%     Unlike time-dependent evolution problems, optimization problems do not require entire approximate solution trajectories with high accuracy, and instead one seeks only some final time point with a small gradient. 
%     %   
%     Thus, in this paper we introduce a restarted MGRIT (R-MGRIT) algorithm, which runs a small number of MGRIT iterations on an approximate trajectory before restarting and moving on to a new trajectory. Each new trajectory corresponds to solving a new all-at-once system starting from an initial point that is the final time point of the previous trajectory.
%     %
%     We find that trajectories can be shortened without increasing the number of trajectories needed for convergence, which increases the estimated parallel speedup while reducing the number of processors required. We provide theoretical results on the convergence of R-MGRIT for linear problems and show that this convergence behavior holds for both linear and nonlinear problems in our numerical experiments.

%    \ok{still finding the abstract not totally clear.}
% \end{abstract}

\begin{abstract}
Previously, we developed a ``parallel-in-iteration'' optimization framework that uses multigrid reduction-in-time (MGRIT) to parallelize slowly converging optimization methods, such as gradient descent. Unlike time-dependent evolution problems, optimization requires only a final point with a sufficiently small gradient, rather than an accurate approximation of the entire trajectory of iterates. This motivates a different use of MGRIT: instead of repeatedly improving the approximation of a single long trajectory, we introduce restarted MGRIT (R-MGRIT), which performs one MGRIT iteration on a trajectory, takes the approximate final point as a new initial guess, and repeats the process on a new trajectory. Each restart thus defines a new all-at-once system, allowing progress toward the solution without accurately approximating the preceding trajectory. Moreover, we find that these trajectories can be shortened without increasing the number of restarts needed for convergence. We provide theoretical results comparing MGRIT and R-MGRIT for linear problems and demonstrate the convergence behavior of the shortened trajectories numerically for both linear and nonlinear problems. Shortening the trajectories increases the estimated parallel speedup while reducing the number of processors required.
\end{abstract}

\section{Introduction}\label{sec:introduction}

Consider the unconstrained optimization problem of minimizing a convex function $F:\mathbb{R}^N \to \mathbb{R}$, i.e., solving $\min_{\mathbf{u} \in \mathbb{R}^N} F(\mathbf{u})$.
Suppose that for solving this problem we employ a 1-step iterative optimization algorithm of the form
\begin{equation}\label{eq:iter_method}
    \mathbf{u}_{i}=\Phi(\mathbf{u}_{i-1})+\mathbf{w}_{i}, \quad i=0,1,\ldots,
\end{equation}
    where $\Phi(\cdot)$ is a linear or nonlinear operator and $\mathbf{w}_{i}$ is a generic forcing term that does not depend on $\mathbf{u}_i$. For example, if $F\coloneqq f+g$, where $f$ is $L$-smooth and $g$ is nonsmooth, then we may use proximal gradient \cite{Beck2009,Beck2017} as the optimization method, which takes the form of \cref{eq:iter_method} with 
    \begin{equation}\label{eq:phi_pg}
     \Phi(\cdot)=(I+s \partial g)^{-1}(I-s\nabla f)(\cdot),
    \end{equation}
    where $\partial g$ denotes the subgradient of $g$, and $\mathbf{w}_i=\mathbf{0}$ for all $i>0$. If $F\coloneqq f$, i.e., the objective function is $L$-smooth, then \cref{eq:phi_pg} reduces to the gradient descent operator \cite{Nesterov1983,Nesterov2004}, another prototypical method of the form \cref{eq:iter_method}. 
These algorithms may converge slowly; for general convex problems, gradient descent and proximal gradient have a convergence rate of $\mathcal{O}(1/i)$, such that tens of thousands or more iterations \cite{Chambolle2016, ODonoghue2013, Beck2009} may be required to get an accurate approximation to the optimum. Since the algorithms are sequential in $i$, a long wall-clock time may then be required to solve the problem.
    
    Assume $N_t$ iterations of the method \cref{eq:iter_method} starting from an initial guess $\mathbf{u}_0 \in \mathbb{R}^{N}$ are required in order to approximate $\mathbf{u}^*$, a minimizer of $F$, to the desired accuracy. Computing this sequence of iterates is equivalent to solving the ``all-at-once'' system,
	\begin{equation}\label{eq:Au=g}
		\mathcal{A}(\mathbf{u})=\left[\begin{array}{cccc}
			I & & & \\
			-\Phi(\cdot) & I & & \\
			& \ddots & \ddots & \\
			& & -\Phi(\cdot) & I
		\end{array}\right]\left[\begin{array}{c}
			\mathbf{u}_0 \\
			\mathbf{u}_1 \\
			\vdots \\
			\mathbf{u}_{N_t}
		\end{array}\right]=\left[\begin{array}{c}
			\mathbf{w}_0 \\
			\mathbf{w}_1 \\
			\vdots \\
			\mathbf{w}_{N_t}
		\end{array}\right]=\mathbf{w},
	\end{equation}
    where $\mathbf{w}_0 \equiv \mathbf{u}_0$ is the initial guess. Parallel-in-time algorithms, such as multigrid reduction-in-time (MGRIT) \cite{Falgout2014}, are designed to solve systems of the form \cref{eq:Au=g} for time-dependent evolution equations \cite{Gander2015, Ong2020}.
    Further inspired by the fact that certain optimization algorithms can be interpreted as discretizations of time-dependent ODEs \cite{Su2016,Attouch2016, Romero2022}, in \cite{Araujo2026} we re-purposed parallel-in-time methods, namely MGRIT, to parallelize over iterations of optimization algorithms by writing those iterations in the form of \cref{eq:Au=g}. We showed that this approach can substantially reduce the predicted wall-clock time compared with executing the optimization algorithm sequentially.
    %We showed that our ``parallel-in-iteration'' approach solves all-at-once systems of the form \cref{eq:Au=g} that correspond to discrete optimization problems with substantially reduced predicted wall-clock time required to solve them. 
    We note that related ideas have also since been used to parallelize gradient flows with Parareal \cite{sun2026}. 

    While parallel-in-time methods approximate entire exact trajectories with high accuracy, the goal for optimization problems is to sufficiently reduce the gradient of the initial guess in order to obtain an approximate minimizer of the objective function. In our parallel-in-iteration approach \cite{Araujo2026}, MGRIT is used to approximate the exact sequence of iterates produced by a fixed number of steps of the sequential optimization method \cref{eq:iter_method}. 
    Since the intermediate iterates generated by \cref{eq:iter_method} are irrelevant to the final solution, expending computational resources to accurately approximate the entire trajectory is inefficient. Instead, we can consider multiple trajectories, and approximate them inexactly with just one or a few MGRIT iterations, where each new trajectory is started from the final point of the previous trajectory.

   We call this restarted MGRIT algorithm R-MGRIT. 
   Given an optimization problem with minimizer $\mathbf{u}^*$ and optimization method \cref{eq:iter_method}, we consider applying MGRIT to the associated all-at-once system $\mathcal{A}(\mathbf{u})=\mathbf{w}$ \cref{eq:Au=g} with $N_t+1$ time points. After one MGRIT iteration, we obtain approximate solutions $\mathbf{v}_i \approx \mathbf{u}_i$, $i=0,1,\ldots,N_t$. If the approximation at the final point, $\mathbf{v}_{N_t}$, is not sufficiently accurate, we formulate a new all-at-once system,
    \begin{equation}\label{eq:upd_Au=g}
        \mathcal{A}(\widetilde{\mathbf{u}})=\widetilde{\mathbf{w}},
    \end{equation}
    where the new initial guess is updated to $\widetilde{\mathbf{w}}_0 = \widetilde{\mathbf{u}}_0 \coloneqq \mathbf{v}_{N_t}$. This updated system determines a new trajectory yielded by \cref{eq:iter_method} starting from the improved initial guess $\widetilde{\mathbf{u}}_0$.
    Rather than accurately approximating a single long trajectory, R-MGRIT approximates a sequence of potentially shorter trajectories, using the final point of each to initialize the next. We therefore refer to R-MGRIT as a \textit{multi-trajectory} parallel-in-iteration algorithm. 
    
    Our motivation is to investigate whether doing a one MGRIT iteration per trajectory leads to a more efficient parallel-in-iteration method than running MGRIT to convergence on a single trajectory, assuming the trajectory has enough optimization steps to reach the desired optimization accuracy.  
    We prove, based on existing MGRIT error propagation theory, that standard MGRIT and R-MGRIT yield approximately the same gradient reduction after a fixed number of iterations. However, we find that we can shorten the consecutive trajectories without increasing the number of trajectories needed for convergence, which increases the parallel speedup while at the same time reducing the number of processors required. 

    The rest of this paper is organized as follows. 
    In \cref{sec:model_problems}, we introduce our model problems and the optimization methods considered in this paper. In \cref{sec:mgrit}, we present the MGRIT algorithm in an optimization context and define our standard and multi-trajectory parallel-in-iteration schemes. In \cref{sec:asympt_equiv}, we compare the standard and restarting MGRIT schemes theoretically, which justifies our approach of truncating the time domain for multi-trajectory solves. \Cref{sec:num_res} provides numerical results, in which we observe that both standard and restarted schemes converge in approximately the same number of iterations; moreover, the number of R-MGRIT trajectories required for convergence is either preserved or reduced by truncating the time domain to a fixed fraction of the original number of time points. In \cref{sec:speedup}, we provide theoretical speedup estimates for the truncated restarting scheme, showcasing its superiority over the standard MGRIT scheme. In \cref{sec:anderson}, we apply Anderson acceleration to accelerate R-MGRIT, and concluding remarks are given in \Cref{sec:conclusions}.

    \section{Model problems}\label{sec:model_problems}
In the following, we briefly summarize our model problems and present the optimization algorithms of our choice for each problem; for a more complete discussion, see \cite{Araujo2026}.
The optimization algorithms are employed in the parallel-in-iteration framework, such that in what follows we speak of a ``fine-grid'' optimization method, and a complementary ``coarse-grid'' optimization method; see discussion in \cref{sec:pii-mgrit}.
    We let $\mathbf{x} \in \mathbb{R}^d$ denote a coordinate on $\Omega \subset \mathbb{R}^d$, such that $\mathbf{x} = x$ when $d = 1$ and $\mathbf{x} = (x, y)$ when $d = 2$.
    The matrices which appear in the definition of our model problems arise from finite-difference discretizations of (negative) Laplacian operators, and are defined as 
    \begin{subequations}
    \label{eq:laplacian}
    \begin{align}
    \label{eq:1d_laplacian}
        A_1 &= -
			\frac{1}{h^2}
            \begin{bmatrix}
                -2 & 1  \\
				1 & -2 & 1 \\
				& \ddots & \ddots & \ddots \\
				& & 1 & -2
            \end{bmatrix}
             \in \mathbb{R}^{n \times n}, 
             \quad \textrm{and} \quad
             \\
             \label{eq:2d_laplacian}
             A_2 &= -\frac{1}{h^2}\begin{bmatrix}
				-T & I &  & \\
				I & -T & I &    \\
                & \ddots& \ddots & \ddots \\
				& &   I & -T
			\end{bmatrix} \in \mathbb{R}^{n^2 \times n^2} 
            \end{align}
            \end{subequations}
            for $d = 1$ and $d = 2$, respectively, with $h=1/(n+1)$.
            Here, $T \in \mathbb{R}^{n \times n}$ is a tridiagonal matrix with elements $-1$, $4$ and $-1$ on the sub-, main and super-diagonals, respectively, and $I \in \mathbb{R}^{n \times n}$ is the identity matrix. The 2-norm of $A_1$ is $\left\|A_1\right\|_2=\frac{1}{h^2} 4 \sin ^2\left(\frac{n \pi}{2(n+1)}\right) \approx \frac{4}{h^2}$ for $h \ll 1$, while the 2-norm of $A_2$ is $\left\|A_2\right\|_2=\frac{1}{h^2} 8 \sin ^2\left(\frac{n \pi}{2(n+1)}\right) \approx \frac{8}{h^2}$ for $h \ll 1$. Throughout the paper, we refer to $A_1$ and $A_2$ as 1D and 2D Laplacian matrices, respectively. 

	\subsection{Quadratic minimization problem (MP1)}\label{subsec:mp1}

Our first model problem, MP1 \cref{eq:mp1}, is the following straightforward quadratic minimization,
	\begin{equation}
        \label{eq:mp1}
		\min_{\mathbf{u} \in \mathbb{R}^{N}} \quad \frac{1}{2} \langle A_1\mathbf{u}, \mathbf{u}  \rangle - \langle \mathbf{b}, \mathbf{u} \rangle, 
	\end{equation}
	where $A_1$ is the 1D Laplacian matrix \eqref{eq:1d_laplacian}, such that $N = n$, and $\mathbf{b} \in \mathbb{R}^{N}$ is some prescribed vector.
    We consider solving this problem with a fine-grid gradient descent iteration, which uses the update 
    	\begin{equation}\label{eq:gd_mp1}
	\begin{aligned}
		\mathbf{u}_{i}=\mathcal{G}_{sf}(\mathbf{u}_{i-1})&\coloneqq (I-s\nabla f)(\mathbf{u}_{i-1})
        % \\&=\mathbf{u}_{i-1}-s(A\mathbf{u}_{i-1}-\mathbf{b})
        \\&=(I-sA_1)\mathbf{u}_{i-1}+s\mathbf{b}, \quad i=1,2,\ldots,
	\end{aligned}
	\end{equation}
    with stepsize $s=1/L$ for $L=\|A_1\|_2$.
    On the coarse-grid we are required to step over the $2/L$ stability limit of gradient descent, and as such we use the proximal point method, with updates given by
    \begin{equation}\label{eq:pp_mp1}
    \begin{aligned}
        \mathbf{u}_i = \mathcal{P}_{sf}(\mathbf{u}_{i-1})&\coloneqq(I+s\nabla f)^{-1}(\mathbf{u}_{i-1}) \\
        &= (I+sA_1)^{-1}(\mathbf{u}_{i-1}+s\mathbf{b}), \quad i=1,2,\ldots.
    \end{aligned}
    \end{equation}

	\subsection{Elastic obstacle problem (MP2)}\label{subsec:mp2}

    Our second model problem, MP2 \cref{eq:mp2}, is an approximate linearized discretization of the elastic obstacle problem (EOP) \cite{Caffarelli1998,Rodrigues2014} with an $\ell^1$ penalty \cite{Tran2015}. Given some discretized obstacle $\boldsymbol{\phi} \in \mathbb{R}^N$ and dimension $d \in \{1,2\}$, letting $\mathbf{b}=-A_d\boldsymbol{\phi}$, MP2 is given by
    \begin{equation}
		\underset{\mathbf{u} \in \mathbb{R}^N}{\min} \quad \frac{1}{2}\langle A_d \mathbf{u}, \mathbf{u}\rangle-\langle \mathbf{b}, \mathbf{u}\rangle+\lambda\left\|(-\mathbf{u})_{+}\right\|_1, \label{eq:mp2}
	\end{equation}
	where $A_d$ is the $d$-dimensional Laplacian matrix from \eqref{eq:laplacian}, 
    %(\cref{eq:1d_laplacian} for $d=1$ with $N=n$, \cref{eq:2d_laplacian} for $d=2$ with $N=n^2$)
    $\lambda>0$ is a predefined penalty parameter, and $\sum_i\left(-(\mathbf{u})_i\right){_{+}}=\left\|(-\mathbf{u})_{+}\right\|_1$, where $(a)_{+}=\max \{0, a\}$ and $\|\cdot\|_1$ denotes the $\ell^1$ norm. For sufficiently large $\lambda$, the solution of \cref{eq:mp2} is an approximate solution of the continuous EOP up to a change of variables \cite{Mangasarian1985, Friedlander2008}. 
    
    As in \cite{Araujo2026}, we consider proximal gradient descent as the fine-grid method for \cref{eq:mp2}, with updates taking the form
    \begin{equation}\label{eq:pg_mp2}
        \begin{aligned}
            \mathbf{u}_i= \mathcal{P}_{sg}\circ \mathcal{G}_{sf}(\mathbf{u}_{i-1})&\coloneqq(I+s \partial g)^{-1}(I-s\nabla f)(\mathbf{u}_i) \\
            &=(I+s \partial g)^{-1}\left( (I-sA_d)\mathbf{u}_{i-1}+s\mathbf{b}\right), \quad i=1,2,\ldots,
        \end{aligned}
    \end{equation}
 where $\nabla f$ denotes the gradient of $f\coloneqq \frac{1}{2}\langle A_d \mathbf{u}, \mathbf{u}\rangle-\langle \mathbf{b}, \mathbf{u}\rangle$ and $\partial g$ denotes the subgradient of $g \coloneqq \lambda \|(-\cdot)_+\|_1$; the proximal operator $P_{sg}=(I+s\partial g)^{-1}$ is given component-wise by
    \begin{equation}\label{eq:prox_mp2}
    \left(\mathcal{P}_{s\lambda\left\|(-\cdot)_{+}\right\|_1}(\mathbf{u})\right)_j= \begin{cases}(\mathbf{u})_j+s\lambda, & (\mathbf{u})_j<-s\lambda, \\
    0, & -s\lambda \leq (\mathbf{u})_j \leq 0 \\
    (\mathbf{u})_j, & (\mathbf{u})_j>0\end{cases}.
        \end{equation}
We take $s=1/L$ with $L=\|A_d\|_2$, $d=1,2$.
      Analogously to MP1, on coarse grids we must step over the $2/L$ stability limit of the proximal gradient descent method, so that on coarse grids we use the alternating proximals method:
        \begin{equation}\label{eq:ap_mp2}
        \begin{aligned}
            \mathbf{u}_i\coloneqq \mathcal{P}_{sg}\circ \mathcal{P}_{sf}(\mathbf{u}_{i-1})&=(I+s \partial g)^{-1}\circ(I+s\nabla f)^{-1}(\mathbf{u}_i) \\
            &=(I+s \partial g)^{-1}\left((I+sA_d-s\mathbf{b})^{-1}(\mathbf{u}_i)\right), \quad i=1,2,\ldots.
        \end{aligned}
    \end{equation}

    \section{Parallel-in-iteration optimization methods}\label{sec:mgrit}
    
    In \cite{Araujo2026}, we presented our parallel-in-iteration framework for optimization problems based on re-purposing parallel time integration algorithms for time-dependent differential equations, motivated by the fact that optimization algorithms often have interpretations as discretizations of time-dependent differential equations (such as gradient flow). In this section, we recall the MGRIT algorithm \cite{Falgout2014}, which we present in the context of optimization. Then, we define our multi-trajectory R-MGRIT parallel-in-iteration scheme. While we use MGRIT in our exposition of the multi-trajectory framework, we note that other parallel-in-time methods, such as PFASST \cite{Minion2010, Emmett2012}, for example, could be employed instead.
    
    \subsection{The MGRIT algorithm} \label{sec:pii-mgrit}
    
\begin{figure}[h!]
            \begin{center}
		\includegraphics[width=0.9\textwidth, trim=0cm 0.1cm 0cm 0.1cm, clip]{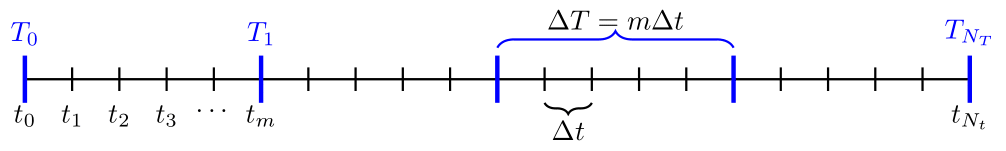}
		\caption{Fine- and coarse-grid discretization meshes with coarsening factor $m$. C-points (in blue) are denoted by $\{ T_j \}$, while all remaining points are F-points (in black); the fine grid is the union of C- and F-points.}
        %are points $t_i$ for $i=1,\ldots,N_t$ such that $i \not = jm$ for $j=1,\ldots,N_T=\frac{N_t}{m}$, and C-points (in blue) are points $T_j=t_{jm}$ for $j=1,\ldots,N_T=\frac{N_t}{m}$.}
        %
		\label{fig:fine_coarse_grids}
            \end{center}
	\end{figure}

    Given an initial guess $\mathbf{u}_0 \in \mathbb{R}^N$ and 1-step iterative method of the form \cref{eq:iter_method}, we consider the associated all-at-once system of the form \cref{eq:Au=g}, where $\mathbf{w}_0=\mathbf{u}_0$ is the initial guess and $\mathbf{u}$ is the exact solution.
    The associated residual vector, given an approximation $\mathbf{v} \approx \mathbf{u}$, is defined as
    \begin{align}\label{eq:mgrit_res}
        \mathbf{r}(\mathbf{v})
        =
        \mathbf{w}
        -
        {\cal A}( \mathbf{v} )=
        \begin{bmatrix}
            \mathbf{w}_0 - \mathbf{v}_0 \\
            \mathbf{w}_1 - [ \mathbf{v}_1 - \Phi( \mathbf{v}_0 )] \\
            \vdots \\
            \mathbf{w}_{N_t} - [ \mathbf{v}_{N_t} - \Phi( \mathbf{v}_{N_t-1} )]
        \end{bmatrix}.
    \end{align}
    We assume that $N_t$ is sufficiently large such that $\mathbf{u}_{N_t}$ is a sufficiently accurate approximation of the minimizer $\mathbf{u}^*$.
	Instead of solving \cref{eq:Au=g} sequentially with time-stepping or forward substitution, MGRIT aims to solve it iteratively and in parallel via the introduction of a coarser approximation of the problem.
Let $t_i=i \Delta t$, $i= 0,1, \ldots, N_t$ be a fine temporal mesh with constant spacing $\Delta t$; in an optimization setting, we take it to be the fixed stepsize
$s>0$ of the underlying optimization method \cref{eq:iter_method}. Then, given an integer coarsening factor $m>1$, we introduce a coarse temporal mesh $T_j=j \Delta T$, $j=0,1, \ldots, N_T$, where $N_T=N_t/m$, with constant spacing $\Delta T=m \Delta t$. We denote the $N_T+1$ points appearing exclusively on the coarse mesh by \emph{C-points}, while all remaining points are \emph{F-points}. 
    \Cref{fig:fine_coarse_grids} illustrates the fine- and coarse-grid temporal meshes.
	
    %
    % An F-relaxation is the process where, starting from a C-point, the time-stepping operator is used to update the values of the $m-1$ F-points that follow it (that is, until the F-point that immediately precedes the next C-point). 
    % %
    % C-relaxation is the update of each C-point by stepping from the preceding F-point using $\Phi$.
    % %
    % FCF-relaxation is simply the application of F-,C- and F-relaxations in this order. 
    % %
    % The post relaxation is taken to be an F-relaxation.
%
    
    MGRIT has both linear \cite{Falgout2014} and Full Approximation Scheme (FAS) variants \cite{Howse2019}, suitable for linear and nonlinear $\Phi$, respectively, and these variants differ in the way that coarse-grid correction is performed. 
    First, we consider the linear case, in which we assume $\Phi$ is a linear operator, such that $\Phi(\mathbf{u}_{i-1})=\Phi\mathbf{u}_{i-1}$ is a matrix-vector product with matrix $\Phi$. We also write the action of the space-time operator $\mathcal{A}(\cdot)$, meaning it operates simultaneously on all spatial coordinates and time points of the solution $\mathbf{u}$, as a matrix-vector product ${\cal A}(\mathbf{u}) = {\cal A}\mathbf{u}$.
	From \cref{eq:Au=g}, we obtain 
	an equivalent system written only in terms of the C-points, which we write in block-matrix form as
    \begin{equation}\label{eq:A_Delta}
    \mathcal{A}_{\Delta} \mathbf{u}_{\Delta}=\left[\begin{array}{cccc}
			I & & & \\
			-\Phi^m & I & & \\
			& \ddots & \ddots & \\
			& & -\Phi^m & I
		\end{array}\right]\left[\begin{array}{c}
			\mathbf{u}_{0} \\
			\mathbf{u}_{m} \\
			\vdots \\
			\mathbf{u}_{N_T m}
		\end{array}\right]=\left[\begin{array}{c}
			\mathbf{w}_0 \\
			\widetilde{\mathbf{w}}_m \\
			\vdots \\
			\widetilde{\mathbf{w}}_{N_Tm}
		\end{array}\right]=\mathbf{w}_\Delta,
    \end{equation}
    where $\widetilde{\mathbf{w}}_{j m}=\mathbf{w}_{j m}+\Phi \mathbf{w}_{j m-1}+\cdots+\Phi^{m-1} \mathbf{w}_{(j-1) m+1}$; note that $\Phi^m$ plays the role of time-stepping operator for this coarse system.
	While system \eqref{eq:A_Delta} has a factor of $m$ fewer equations than \eqref{eq:Au=g}, it is no less expensive to solve because (naively) computing the action of $\Phi^m$ is $m$ times as expensive as computing that of $\Phi$.
    So, instead, MGRIT solves an approximation to the \emph{ideal} coarse-grid system \eqref{eq:A_Delta} arising from  replacing the \emph{ideal} coarse-grid time-stepping operator $\Phi^m$ with an approximation:
	\begin{equation} \label{eq:phi_delta}
			\Psi \approx \Phi^m. 
	\end{equation}
    Then, MGRIT solves the approximate system $\mathcal{B}_{\Delta} \mathbf{u}_{\Delta}=\mathbf{w}_{\Delta}$, where the approximation $\mathcal{B}_{\Delta} \approx \mathcal{A}_{\Delta}$ replaces the $\Phi^m$ terms in \eqref{eq:A_Delta} by $\Psi$.
   
    Throughout the paper, we assume MGRIT is initialized with an initial approximate solution vector of the form $\mathbf{v}^{(0)}=\begin{bmatrix}
        \mathbf{u}_0 & \mathbf{u}_0 & \cdots & \mathbf{u}_0
    \end{bmatrix}^T$ which simply extends the initial guess $\mathbf{u}_0$ to all other time points, which will be important to maintain progress in the restarted method. Letting $\mathbf{r}^{(0)}$ denote the initial MGRIT residual \cref{eq:mgrit_res}, we define the linear MGRIT algorithm in \cref{alg:mgrit} with the standard MGRIT pre- and post-relaxation schemes \cite{Falgout2014,Dobrev2017}, where $R_I(\cdot)$ stands for the injection operator. Note that we omit the post F-relaxation that follows the correction at C-points, since it is made redundant by the FCF-relaxation in the following iteration.
    \begin{algorithm}
    \caption{2-level MGRIT($\Phi, \Psi, \mathbf{u_0}, \mathbf{w}$,$\varepsilon$).}\label{alg:mgrit}
    
    \begin{algorithmic}[1]
        \State Initialize $\mathbf{u}$ given initial guess $\mathbf{u}_0 \in \mathbb{R}^N$: $\mathbf{u}=\begin{bmatrix} \mathbf{u}_0 & \cdots & \mathbf{u}_0 \end{bmatrix}^T \in \mathbb{R}^{N(N_t+1)}.$
    \Repeat
        \State Apply FCF-relaxation to $\mathcal{A}\mathbf{u}=\mathbf{w}$ using $\Phi$.
        \State Inject residual to the coarse grid: $\mathbf{r}_\Delta=R_I(\mathbf{w}-\mathcal{A}\mathbf{u})$.
        
        \State Solve coarse-grid correction problem using $\Psi$: $\mathcal{B}_\Delta\mathbf{e}_\Delta=\mathbf{r}_\Delta$.
        \State Correct approximate solution at C-points: $\mathbf{u}_\Delta \leftarrow \mathbf{u}_\Delta + \mathbf{e}_\Delta$.
    \Until{ $\|\mathbf{r}\|/\|\mathbf{r}^{(0)}\|\leq \varepsilon$.}
    \end{algorithmic}
    \end{algorithm}
    We note that while the MGRIT convergence criterion is given in terms of the residual \cref{eq:mgrit_res}, in \cite[Lemma 5.1]{Araujo2026} we proved that reducing the norm of the initial residual by a factor $\varepsilon$ implies that the norm of the gradient is also reduced to a similar factor; i.e., $\|\nabla f(\mathbf{v}_{N_t})\|/\|\nabla f(\mathbf{u}_{0})\|\approx \varepsilon$, where $\mathbf{v}_{N_t}$ denotes the MGRIT approximate solution at the final time point $N_t$.
   
    Now let us consider the MGRIT coarse-grid correction in the more general FAS case, where $\Phi(\cdot)$ in \eqref{eq:Au=g} is not (necessarily) a linear operator. 
    Let $\mathbf{v} \approx \mathbf{u}$ be an approximate solution of $\mathcal{A}(\mathbf{u})=\mathbf{w}$ in \cref{eq:Au=g}.   
    After performing the FCF-relaxation step, both the fine-grid approximation $\mathbf{v}$ and nonlinear residual $\mathbf{w}-\mathcal{A}(\mathbf{v})$ are injected to the coarse-grid:
	\begin{equation}
		\mathbf{v}_{\Delta}=R_{I}(\mathbf{v}), \quad \mathbf{r}_{\Delta}=R_{I}(\mathbf{w}-\mathcal{A}(\mathbf{v})).
	\end{equation}
    The \emph{ideal} FAS coarse-grid problem is then to solve the system
    %is a C-point Schur complement \hds{Is there a nonlinear notion of a Schur complement?} of the residual equation \cref{eq:mgrit_res} \hds{this is not clear},
    \begin{equation}\label{eq:c_mgrit_problem_fas}		 \mathcal{A}_{\Delta}\left(\mathbf{u}_{\Delta}\right)=\mathcal{A}_{\Delta}\left(\mathbf{v}_{\Delta}\right)+\mathbf{r}_{\Delta}
    \end{equation}
    %\hds{How does this relate to (4.4)? This is the nonlinear version of (4.10).}
    for the coarse-grid solution $\mathbf{u}_{\Delta}$, where similarly to \eqref{eq:A_Delta},
    \begin{equation} 
		\mathcal{A}_\Delta(\cdot)=\left[\begin{array}{cccc}
			I & & & \\
			-\Phi^m(\cdot) & I & & \\
			& \ddots & \ddots & \\
			& & -\Phi^m(\cdot) & I
		\end{array}\right].
	\end{equation}
    The coarse-grid system \eqref{eq:c_mgrit_problem_fas} results from the elimination of F-points variables from the residual equation ${\cal A}( \mathbf{u} ) - {\cal A}(\mathbf{v}) = \mathbf{r}( \mathbf{v} )$.
    Similarly to the linear case, the MGRIT iteration solves an approximation of \cref{eq:c_mgrit_problem_fas}, 
    \begin{equation} 
    \mathcal{B}_{\Delta}\left(\mathbf{u}_{\Delta}\right) = \mathcal{B}_{\Delta}\left(\mathbf{v}_{\Delta}\right)+\mathbf{r}_{\Delta},
    \end{equation} 
    where $\mathcal{B}_{\Delta}(\cdot) \approx \mathcal{A}_\Delta(\cdot)$ replaces the $\Phi^m(\cdot)$ terms in $\mathcal{A}_\Delta(\cdot)$ by a computationally cheaper $\Psi(\cdot) \approx \Phi^m(\cdot)$. 

    We define the 2-level FAS-MGRIT iteration in \cref{alg:fas-mgrit}, once again omitting the redundant post-F relaxation.
    \begin{algorithm}
    \caption{2-level FAS-MGRIT($\Phi, \Psi, \mathbf{u_0}, \mathbf{w}$,$\varepsilon$).}\label{alg:fas-mgrit}
    \begin{algorithmic}[1]
        \State Initialize $\mathbf{u}$ given initial guess $\mathbf{u}_0 \in \mathbb{R}^N$: $\mathbf{u}=\begin{bmatrix}
        \mathbf{u}_0 & \cdots & \mathbf{u}_0     \end{bmatrix}^T \in \mathbb{R}^{N(N_t+1)}.$
    \Repeat
        \State Apply FCF-relaxation to $\mathcal{A}(\mathbf{u})=\mathbf{w}$ using $\Phi$.
        \State Inject approximate solution and residual to coarse grid: 
        $\mathbf{u}_\Delta=R_I(\mathbf{u})$,  $\mathbf{r}_\Delta=R_I(\mathbf{w}-\mathcal{A}(\mathbf{u}))$.
        \State Solve coarse-grid correction problem: $\mathcal{B}_\Delta(\widehat{\mathbf{u}}_\Delta)=\mathcal{B}_\Delta(\mathbf{u}_\Delta)+\mathbf{r}_\Delta$.
        \State Compute coarse-grid error approximation: $\mathbf{e}_\Delta=\widehat{\mathbf{u}}_\Delta-\mathbf{u}_\Delta$.
        \State Correct approximate solution at C-points: $\mathbf{u}_\Delta \leftarrow \mathbf{u}_\Delta + \mathbf{e}_\Delta$.
    \Until{ $\|\mathbf{r}\|/\|\mathbf{r}^{(0)}\|\leq \varepsilon$.}
    \end{algorithmic}
    \end{algorithm} Similarly to the linear case, in \cite[Lemma 5.1]{Araujo2026} we proved that $\|\mathbf{r}\|/\|\mathbf{r}^{(0)}\|\leq \varepsilon$ implies that  $\|G_{sF}(\mathbf{v}_{N_t})\|/\|G_{sF} (\mathbf{u}_{0})\| \approx \varepsilon$, where $G_{sF}(\cdot)$ denotes the generalized gradient of $F$, given by
	\begin{equation} \label{eq:G_sF}
		G_{sF}\left(\mathbf{v}_i\right):=\frac{1}{s}\left[\mathbf{v}_i-\mathcal{P}_{sg}\left(\mathbf{v}_i-s \nabla f\left(\mathbf{v}_i\right) \right)\right]=\frac{1}{s}(I-\mathcal{P}_{sg} \circ \mathcal{G}_{sf}(\mathbf{v}_i)) 
	\end{equation}
    for some $\mathbf{v}_i \in \mathbb{R}^{N}$.
    Note that $\mathcal{P}_{sg} \circ \mathcal{G}_{sf}$ is the same as the proximal gradient operator \cref{eq:pg_mp2}; in particular, $\mathcal{P}_{sg} \circ \mathcal{G}_{sf}=I-sG_{sF}$.
    
    While we have described 2-level MGRIT methods, multilevel MGRIT schemes are  obtained by recursively applying the 2-level schemes. In our numerical tests, we apply the standard MGRIT scheme for MP1 \cref{eq:mp1} and the FAS-MGRIT scheme for MP2 \cref{eq:mp2}. For both MP1 and MP2, MGRIT is iterated until the 2-norm of the residual \cref{eq:mgrit_res} is reduced by 8 orders of magnitude from its initial value. 

  \subsection{The restarted MGRIT algorithm for optimization problems (R-MGRIT)}

We now describe the R-MGRIT algorithm, \cref{alg:mgrit_rest}. As introduced in \cref{sec:introduction}, instead of approximating a single exact trajectory yielded by the underlying optimization method \cref{eq:iter_method}, we treat the initial guess $\mathbf{u}_0$ of the all-at-once system $\mathcal{A}(\mathbf{u})=\mathbf{w}$ \cref{eq:Au=g} as a parameter to be reinitialized. Each new trajectory is approximated using one MGRIT iteration.\footnote{In general, more than one MGRIT iteration can be performed if needed. However, in our experiments we have found that using more than one provides no additional computational benefit.} Thus, R-MGRIT can be interpreted as a restarting scheme that approximates multiple trajectories, each hopefully starting from a point closer to the true optimum. It is important to propagate $\mathbf{u}_0$ across the C-points in the new trajectory because a random initialization of the C-points may reintroduce large gradients. We note that one R-MGRIT iteration costs the same as one MGRIT iteration from \cref{alg:mgrit,alg:fas-mgrit}.
     % \ok{as a general rule, abbreviations should be defined once initially and used consistently} \ok{so i guess i'm confused, you state here you're going to describe the approach, but all of the above reads as though the reader has retained all of the discussion from the intro? (not a good idea, too much of a gap since then)}
    % \ok{some of the alg does not make sense: 1. what is $A$, do you mean ${\cal A}$? and you pass parameters into mgrit that fully specify the system already, so i don't think you can say ``run MGRIT<blah parameters> on <blah system>'' because that's overloading the notation; also, doesnt passing a tol epsilon mean iterate mgrit until that tol, rather than run 1 iter?}
      \begin{algorithm}
    \caption{R-MGRIT($\Phi, \Psi, \mathbf{u}_0, \mathbf{w}$,$\varepsilon$,$N_t$).}\label{alg:mgrit_rest}
    \begin{algorithmic}[1]
    \Repeat
        \State Initialize $\mathbf{u}$ given initial guess $\mathbf{u}_0 \in \mathbb{R}^N$: $\mathbf{u}=\begin{bmatrix}
        \mathbf{u}_0 & \cdots & \mathbf{u}_0     \end{bmatrix}^T \in \mathbb{R}^{N(N_t+1)}.$
        \State Run 1 iteration of MGRIT on all-at-once system $\mathcal{A}(\mathbf{u})=\mathbf{w}$.
        \State Set $\mathbf{u}_0 \longleftarrow \mathbf{u}_{N_t}$, $\mathbf{w}_0 \longleftarrow \mathbf{u}_0$.
    \Until{ $\|G_{sF}(\mathbf{u}_{N_t})\|/\|G_{sF}(\mathbf{u}_{0})\| \leq \varepsilon$.}
    \end{algorithmic}
    \end{algorithm}

\section{Asymptotic final-point analysis of MGRIT and R-MGRIT}\label{sec:asympt_equiv}

In the following, we compare standard linear MGRIT and R-MGRIT theoretically. We use the index $k$ to denote the number of iterations for standard MGRIT and the number of restarts (i.e., trajectories) for R-MGRIT. 
Based on linear MGRIT error propagation theory \cite{Dobrev2017,Southworth2019}, we demonstrate that if the approximate solutions at the final time point $N_t=N_Tm$ of the $k$th MGRIT iteration and $k$th R-MGRIT trajectory are approximately the same, converging to the same asymptotic limit as $N_T \to \infty$. 

Consider an iterative method of the form \cref{eq:iter_method} with a linear operator $\Phi$. Furthermore, we assume $\mathbf{w}_i=\mathbf{c}$ for $i>0$, such that, given an initial guess $\mathbf{u}_0$, we have
\begin{equation}
            \mathbf{u}_i=\Phi\mathbf{u}_{i-1}+\mathbf{c}, \quad i=1,2,\ldots
\end{equation}
Writing the $i$th iterate yielded by \cref{eq:iter_method} in terms of the initial guess $\mathbf{u}_0$ yields
\begin{equation}
    \begin{aligned}
        % &=\Phi(\Phi\mathbf{u}_{j-2}+\mathbf{c})+\mathbf{c} \\
        % &=\Phi^2\mathbf{u}_{j-2}+(I+\Phi)\mathbf{c} \\
        \mathbf{u}_i=\Phi^{i}\mathbf{u}_0 +(I+\Phi+\cdots+\Phi^{i-1})\mathbf{c}, \quad i=1,2,\ldots
    \end{aligned}
\end{equation}
%Without loss of generality, by shifting the coordinate system to $\mathbf{u} - \mathbf{u}^*$, we may assume $\mathbf{c}=\mathbf{0}$ and $\mathbf{u}^*=\mathbf{0}$. 
Since $\mathbf{u}^*$ is a fixed point of the iteration,
i.e., $\mathbf{u}^*=\Phi\mathbf{u}^*+\mathbf{c}$,
we may shift coordinates by defining
$\mathbf{v}_j=\mathbf{u}_j-\mathbf{u}^*$, with the resulting iteration being
$\mathbf{v}_i=\Phi\mathbf{v}_{i-1}$.
Thus, without loss of generality, we simplify the above iteration to 
\begin{equation}\label{eq:homog_gd}
\mathbf{u}_i = \Phi\mathbf{u}_{i-1}=\Phi^i\mathbf{u}_0, \quad i=1,2,\ldots
\end{equation}
For concreteness, we may consider $\Phi\coloneqq I-sA_1$ in the case of gradient descent \cref{eq:gd_mp1} for MP1, \cref{eq:mp1}. 

In the following, we first explicitly determine the MGRIT error propagation matrix for $k$ MGRIT iterations. Assuming a constant initialization strategy, we then characterize the MGRIT and R-MGRIT approximate solutions, for the case of $k$ MGRIT iterations and $k$ R-MGRIT trajectories with 1 MGRIT iteration per trajectory, at the final time point in terms of the MGRIT parameters. Finally, we prove that the MGRIT and R-MGRIT end-point approximate solutions converge to the same asymptotic limit for fixed $k$ as $N_T \to \infty$. 

\subsection{MGRIT error propagation after $k$ iterations}
Let $\Phi$ and $\Psi$ denote the fine- and coarse-grid operators. After $k$ iterations, the C-point error propagation block matrix $\mathcal{E}$ propagates the error at C-points via \begin{equation}
    \mathbf{e}_{c}^{(k)}=\mathcal{E}^k\mathbf{e}_{c}^{(0)},
\end{equation}
where $\mathbf{e}_{c}^{(0)}$ and $\mathbf{e}_{c}^{(k)}$ denote the C-point initial error and the C-point error after $k$ iterations, respectively. 
For MGRIT with FCF-relaxation, $\mathcal{E}$ is given by \cite{Dobrev2017,Southworth2019}
\begin{equation}\label{eq:error_prop_matrix_def}
    \mathcal{E} = \begin{bmatrix}
    \mathbf{0} & \mathbf{0} & \cdots & \mathbf{0} &  \mathbf{0}\\
    \mathbf{0} & \mathbf{0} & \cdots & \mathbf{0} &  \mathbf{0}\\
    (\Phi^m-\Psi)\Phi^m & \mathbf{0} & \cdots & \mathbf{0} &  \mathbf{0} \\
    \Psi(\Phi^m-\Psi)\Phi^m & (\Phi^m-\Psi)\Phi^m & \cdots & \mathbf{0} &  \mathbf{0}\\
    \vdots & \vdots & \ddots & \vdots & \vdots\\
    \Psi^{N_T-2}(\Phi^m-\Psi)\Phi^m & \Psi^{N_T-3}(\Phi^m-\Psi)\Phi^m & \cdots & \mathbf{0} &  \mathbf{0}
    \end{bmatrix},
\end{equation}
or block-wise, by
\begin{equation}
    (\mathcal{E})_{i,j} = \begin{cases} 
   \Psi^{i-j-2} (\Phi^m-\Psi)\Phi^m & \text{if } i-j \ge  2, \\ 
    \mathbf{0} & \text{otherwise.}
    \end{cases}
\end{equation}
By taking the $k$th power of $\mathcal{E}$ in \cref{eq:error_prop_matrix_def}, we explicitly determine the error propagation matrix $\mathcal{E}^k$ for $k$ MGRIT iterations.
We now introduce an assumption to hold uniformly throughout this section; while not strictly necessary for all of the results, applying the assumption uniformly aids in the presentation.
\begin{assumption}\label[assumption]{ass:operators}
The fine- and coarse-grid operators $\Phi$ and $\Psi$ are
simultaneously unitarily diagonalizable and satisfy
$\rho(\Phi)<1$ and $\rho(\Psi)<1$.
\end{assumption}

\begin{proposition}\label{prop:error_prop_matrix_powers}
Let \cref{ass:operators} hold.
For any iteration $k \ge 1$, the block-components of $\mathcal{E}^k$ are given by:
\begin{equation}\label{eq:error_prop_matrix_powers}
    (\mathcal{E}^k)_{i,j} = \begin{cases} 
    \binom{i-j-k-1}{k-1} \Psi^{i-j-2k} ((\Phi^m-\Psi)\Phi^m)^k & \text{if } i-j \ge  2k, \\ 
    0 & \text{otherwise.}
    \end{cases}
\end{equation}
\end{proposition}
\begin{proof}
By induction on $k$. The base case $k=1$ trivially matches \cref{eq:error_prop_matrix_def} since $\binom{i-j-2}{0} = 1$. 
Now, assume the formula holds for $k-1$. We compute the $(i,j)$-th block of $\mathcal{E}^k = {\cal E} {\cal E}^{k-1}$ via the block matrix product 
\begin{equation}\label{eq:induc_hyp}
    (\mathcal{E}^k)_{i,j} = \sum_{p=0}^{N_T} (\mathcal{E})_{i,p} (\mathcal{E}^{k-1})_{p,j}.
\end{equation}
By \eqref{eq:error_prop_matrix_powers} holding for ${\cal E}$ and ${\cal E}^{k - 1}$, the addends in this sum are non-zero only if $(\mathcal{E})_{i,p} \neq 0 \iff i -p \ge 2$ and $(\mathcal{E}^{k-1})_{p,j} \neq 0 \iff p -j \ge 2(k-1)$. 
Adding these two inequalities, we have that $(\mathcal{E}^k)_{i,j} = 0$ if $i - j < 2k$.
Otherwise, for $i\ge j+2k$, substituting \cref{eq:error_prop_matrix_powers} into \cref{eq:induc_hyp} yields
\begin{align*}
    (\mathcal{E}^k)_{i,j} 
    &= \sum_{p=j+2k-2}^{i-2} \left( \Psi^{i-p-2} ((\Phi^m-\Psi)\Phi^m) \right) \left({\textstyle \binom{p-j-k}{k-2}} \Psi^{p-j-2k+2} ((\Phi^m-\Psi)\Phi^m)^{k-1} \right)\\
    &=\Psi^{i-j-2k} ((\Phi^m-\Psi)\Phi^m)^k \sum_{p=j+2k-2}^{i-2} \binom{p-j-k}{k-2},   
\end{align*}
where we used the commutativity of $\Phi$ and $\Psi$.
Recall the hockey-stick (or Fermat's) identity for non-negative integers $n \geq r$, \begin{equation}\label{eq:hockey_stick}
    \sum_{i=r}^{n} \binom{i}{r} = \binom{n+1}{r+1}.
\end{equation}
Applying the change of variables $q = p - j - k$, we get that
\begin{equation}
    \sum_{p=j+2k-2}^{i-2} \binom{p-j-k}{k-2}=\sum_{q=k-2}^{i-j-k-2} \binom{q}{k-2}=\binom{i-j-k-1}{k-1},
\end{equation}
where the final equality comes from applying \cref{eq:hockey_stick}.
\end{proof}

From \cref{prop:error_prop_matrix_powers}, we generalize \cite[Lemma 3.2]{Dobrev2017} to determine the C-point errors after $k$ iterations.

\begin{lemma}\label{lem:Cpt_err} 
Let \cref{ass:operators} hold.
After $k$ iterations of MGRIT with FCF-relaxation, the approximate solution at the $i$th C-point satisfies $\mathbf{v}^{(k)}_{im}=\mathbf{u}_{im}-\mathbf{e}^{(k)}_{im}$, where the error at the C-points satisfies
\begin{align}\label{eq:Cpt_err}
    %\mathbf{e}^{(k)}_0&=\mathbf{e}^{(k)}_m=\mathbf{0}, \\
    \mathbf{e}^{(k)}_{im} &= \mathbf{0},
\quad i=0,\ldots,\min(2k-1,N_T), \\
    \mathbf{e}^{(k)}_{im}&=\sum_{j=0}^{i-2k} \binom{i-j-k-1}{k-1} \Psi^{i-j-2k} ((\Phi^m-\Psi)\Phi^m)^k \mathbf{e}^{(0)}_{jm}, \quad i=2k,\ldots,N_T.
\end{align}
\end{lemma}
\begin{proof}
The error at the $i$-th C-point, $\mathbf{e}^{(k)}_{im}$, is determined by the product of the $i$-th block row of $\mathcal{E}^k$ with the initial C-point error vector $\begin{bmatrix}
    \mathbf{e}^{(0)}_0 & \mathbf{e}^{(0)}_m & \cdots & \mathbf{e}^{(0)}_{N_Tm}
\end{bmatrix}^T$:
\begin{equation}\label{eq:ek_jm_sum}
\mathbf{e}^{(k)}_{im} 
% = (\mathcal{E}^k)_{i,0:N_T}\mathbf{e}^{(0)}_c
= \sum_{j=0}^{N_T} (\mathcal{E}^k)_{i,j} \mathbf{e}^{(0)}_{jm}.   
\end{equation}
The statement then follows from substituting \cref{eq:error_prop_matrix_powers} into \cref{eq:ek_jm_sum}.
\end{proof}

\subsection{Characterizing (R-)MGRIT final-point approximate solutions}

We now characterize the approximate solutions at the final time point
$N_t$ for standard MGRIT and R-MGRIT. For both methods, we use the
constant initialization strategy described in \cref{alg:mgrit,alg:mgrit_rest}.

First, consider standard MGRIT applied to a trajectory with initial
guess $\mathbf{u}_0$. Since the exact solution at the $j$th C-point
is $\mathbf{u}_{jm}=\Phi^{jm}\mathbf{u}_0$, the constant initialization
gives an initial error
$\mathbf{e}^{(0)}_{jm}=(\Phi^{jm}-I)\mathbf{u}_0$.
Substituting this into \cref{lem:Cpt_err}, the error after $k$
MGRIT iterations at the final time point is
$\mathbf{e}^{(k)}_{N_t}=E^{(k)}_{N_T}\mathbf{u}_0$, where
\begin{equation}\label{eq:Ek_NT}
 E^{(k)}_{N_T} =
 \sum_{j=0}^{N_T-2k}
 \binom{N_T-j-k-1}{k-1}
 \Psi^{N_T-j-2k}
 ((\Phi^m-\Psi)\Phi^m)^k
 (\Phi^{jm}-I).
\end{equation}
Since the exact solution at the final time point is
$\mathbf{u}_{N_t}=\Phi^{N_t}\mathbf{u}_0$, the corresponding
MGRIT approximation is
\begin{equation}\label{eq:std_vk_Nt_Ek}
    \mathbf{v}^{(k)}_{N_t}
    = \underbrace{(\Phi^{N_t}-E^{(k)}_{N_T})}_{=: V^{(k)}_{N_t}}\mathbf{u}_0.
\end{equation}

Next, consider R-MGRIT with one MGRIT iteration per trajectory.
Since the same constant initialization strategy is used for each
trajectory, the operator mapping its initial guess to its approximate
final point is always the one-iteration MGRIT operator
\begin{equation}\label{eq:restarted_single_iter_op}
 \widetilde{V}_{N_t}
 \coloneqq V^{(1)}_{N_t}
 = \Phi^{N_t}-E^{(1)}_{N_T}.
\end{equation}
Let $\widetilde{\mathbf{u}}^{(k)}_0$ denote the initial guess
of the $k$th R-MGRIT trajectory, with
$\widetilde{\mathbf{u}}^{(1)}_0=\mathbf{u}_0$.
The approximate final point of this trajectory is
$\widetilde{\mathbf{v}}^{(k)}_{N_t}
=\widetilde{V}_{N_t}\widetilde{\mathbf{u}}^{(k)}_0$.
By construction, this becomes the initial guess for the next
trajectory, i.e.,
$\widetilde{\mathbf{u}}^{(k+1)}_0
=\widetilde{\mathbf{v}}^{(k)}_{N_t}$.
Applying this relation recursively, we obtain
\begin{equation}\label{eq:res_vk_Nt_Ek}
 \widetilde{\mathbf{v}}^{(k)}_{N_t}
 =(\widetilde{V}_{N_t})^k\mathbf{u}_0
 =(\Phi^{N_t}-E^{(1)}_{N_T})^k\mathbf{u}_0.
\end{equation}

Thus, standard MGRIT applies $k$ iterations to a single trajectory,
yielding the final-point approximation
$\mathbf{v}^{(k)}_{N_t} = V^{(k)}_{N_t}\mathbf{u}_0$, whereas R-MGRIT applies the same
single-iteration final-point operator $k$ times, yielding
$ \widetilde{\mathbf{v}}^{(k)}_{N_t} = (\widetilde{V}_{N_t})^k\mathbf{u}_0$.
In the following sections, we compare these two approximations
as $N_T\to\infty$.

\subsection{Error operators}\label{subsec:error_oper_sim_diag}

Let $
    Q=\begin{bmatrix}
        \mathbf{q}_1 & \mathbf{q}_2 & \cdots &\mathbf{q}_N
    \end{bmatrix}
    $ be the unitary matrix that diagonalizes both $\Phi$ and $\Psi$ under \cref{ass:operators},
 such that  
 \begin{equation}\label{eq:diag_assumption}
    \Phi=Q \operatorname{diag}(\phi_1,\ldots,\phi_N)Q^*, \quad      \Psi=Q \operatorname{diag}(\psi_1,\ldots,\psi_N)Q^*,
\end{equation}
is satisfied with $\phi_n$ and $\psi_n$ denoting the $n$th eigenvalues of $\Phi$ and $\Psi$, respectively. 
We denote by $u_{0,n}=\mathbf{q}_n^*\mathbf{u}_0$ the $n$th modal component of the initial guess $\mathbf{u}_0$, where $\mathbf{q}_n$ is the $n$th column of $Q$.
Then, after an eigenvector expansion of $\mathbf{e}_{c}^{(1)}=\mathcal{E}\mathbf{e}_{c}^{(0)}$ on the basis of eigenvectors given by the columns of $Q$, it follows that the action of $\mathcal{E}$ on the $n$th mode is given by $\mathcal{E}_n \in \mathbb{C}^{N_T + 1 \times N_T + 1}$, given by \cite{Dobrev2017,Southworth2019}
\begin{equation}\label{eq:error_prop_matrix_mode}
    \mathcal{E}_n = \begin{bmatrix}
    0 & 0 & \cdots & 0 &  0\\
    0 & 0 & \cdots & 0 &  0\\
    (\phi_n^m-\psi_n)\phi_n^m & 0 & \cdots & 0 &  0 \\
    \psi_n(\phi_n^m-\psi_n)\phi_n^m & (\phi_n^m-\psi_n)\phi_n^m & \cdots & 0 &  0\\
    \vdots & \vdots & \ddots & \vdots & \vdots\\
    \psi_n^{N_T-2}(\phi_n^m-\psi_n)\phi_n^m & \psi_n^{N_T-3}(\phi_n^m-\psi_n)\phi_n^m & \cdots & 0 &  0
    \end{bmatrix}.
\end{equation}
Following similar steps, it also follows from \cref{prop:error_prop_matrix_powers} that the action of $E^{(k)}_{N_T}$ from \cref{eq:Ek_NT} on the $n$th mode (i.e., the $n$th modal component of $\mathbf{u}_0 \in \mathbb{R}^N$) is
\begin{equation}\label{eq:Ek_NT_mode}
 E^{(k)}_{N_T,n} = ((\phi_n^m-\psi_n)\phi_n^m)^k\sum_{j=0}^{N_T-2k} \binom{N_T-j-k-1}{k-1} \psi_n^{N_T-j-2k} (\phi_n^{jm}-1).
\end{equation}
such that the modal equivalents of \cref{eq:std_vk_Nt_Ek} and \cref{eq:res_vk_Nt_Ek} are 
\begin{equation}\label{eq:std_vk_Nt_Ek_mode}
  v^{(k)}_{N_t, n} = V^{(k)}_{N_t, n}u_{0,n}=(\phi_n^{N_t} - E^{(k)}_{N_T, n})u_{0,n},  
\end{equation}
and
\begin{equation}\label{eq:res_vk_Nt_Ek_mode}
   \widetilde{v}^{(k)}_{N_t, n} =  \widetilde{V}_{N_t, n}^k u_{0,n}=(\phi_n^{N_t} - E^{(1)}_{N_T, n})^ku_{0,n},
\end{equation}
respectively, where $v^{(k)}_{N_t,n}$ is the $n$th modal component of $\mathbf{v}^{(k)}_{N_t}$, and $\widetilde{v}^{(k)}_{N_t,n}$ is the $n$th modal component of $\widetilde{\mathbf{v}}^{(k)}_{N_t}$. We now establish the asymptotic modal error operator as $N_T \to +\infty$.

\begin{lemma}\label{lem:Ek_infty_mode}
Let \cref{ass:operators} hold.
 As $N_T \to +\infty$, the modal error operator $E^{(k)}_{N_T,n}$ from \cref{eq:Ek_NT_mode} for $k$ MGRIT iterations asymptotically converges to a well-defined limit modal operator
\begin{equation}\label{eq:Ek_infty_mode}
    E^{(k)}_{\infty,n}:=-\left(\frac{(\phi_n^m-\psi_n)\phi_n^m}{1-\psi_n}\right)^k=\lim_{N_T \to +\infty} E^{(k)}_{N_T,n}.
    \end{equation}
\end{lemma}

\begin{proof}
Applying the change of variables $p = N_T - 2k - j$ reverses the summation order, rewriting the modal error operator from \cref{eq:Ek_NT_mode} as:
\begin{equation}\label{eq:Ek_reversed_sum}
    E^{(k)}_{N_T,n} = ((\phi_n^m-\psi_n)\phi_n^m)^k \sum_{p=0}^{N_T-2k} \binom{p+k-1}{k-1} \psi_n^p (\phi_n^{(N_T-2k-p)m}-1).
\end{equation}
Let $M = N_T - 2k$. We can split the summation into two distinct parts and evaluate their limits as $M \to \infty$. For the first part involving the $\phi_n$ term, let
\begin{equation}
    S_M = \sum_{p=0}^{M} \binom{p+k-1}{k-1} \psi_n^p \phi_n^{(M-p)m}.
\end{equation}
To bound this sum, let $\rho = \max(|\psi_n|, |\phi_n^m|)$. By hypothesis, both $|\phi_n|<1$ and $|\psi_n|<1$, which implies $\rho < 1$. Taking the absolute value yields:
\begin{equation}
    |S_M| \le \sum_{p=0}^{M} {\textstyle \binom{p+k-1}{k-1}} |\psi_n|^p |\phi_n^m|^{M-p} \le \sum_{p=0}^{M} {\textstyle \binom{p+k-1}{k-1}} \rho^p \rho^{M-p} = \rho^M {\textstyle \sum_{p=0}^{M} \binom{p+k-1}{k-1}}.
\end{equation}
Using the hockey-stick identity \eqref{eq:hockey_stick}, the sum of the binomial coefficients evaluates to $\binom{M+k}{k}$, giving the bound $|S_M| \le \binom{M+k}{k} \rho^M$.
The expression $\binom{M+k}{k}$ is a polynomial in $M$ of degree $k$, while $\rho^M$ is an exponentially decaying term. Since exponential decay dominates polynomial growth as $M \to \infty$, we have $\lim_{M \to +\infty} \binom{M+k}{k} \rho^M = 0$, meaning $\lim_{M \to +\infty} S_M = 0$.

The second part of the summation in \cref{eq:Ek_reversed_sum},
$-\sum_{p=0}^{M}\binom{p+k-1}{k-1}\psi_n^p$,
is a partial sum of the negative binomial series.
As $M \to \infty$, this yields the Maclaurin series expansion for $-(1-\psi_n)^{-k}$:
\begin{equation}
    \lim_{M \to +\infty} - \sum_{p=0}^{M} \binom{p+k-1}{k-1} \psi_n^p = -(1-\psi_n)^{-k}.
\end{equation}
Because $|\psi_n| < 1$, this series converges absolutely. Substituting these two limits back into the original expression yields
\begin{equation}
    \lim_{N_T \to +\infty} E^{(k)}_{N_T,n} = ((\phi_n^m-\psi_n)\phi_n^m)^k \left( 0 - \frac{1}{(1-\psi_n)^k} \right) = -\left(\frac{(\phi_n^m-\psi_n)\phi_n^m}{1-\psi_n}\right)^k,
\end{equation}
which recovers $E^{(k)}_{\infty,n}$ in \cref{eq:Ek_infty_mode}.
\end{proof}

\subsection{Asymptotic final-point (R-)MGRIT approximate solutions}\label{subsec:asympt_equiv}

We \linebreak now show that for any fixed $k$, both MGRIT with $k$ iterations and R-MGRIT with 1 MGRIT iteration on $k$ consecutive trajectories converge to the same final-point approximate solution, asymptotically for $N_T \to +\infty$.

\begin{theorem}\label{thm:asymptotic_equivalence}
Let \cref{ass:operators} hold. 
For any fixed $k\geq 1$, the MGRIT and R-MGRIT final-point solutions from \cref{eq:std_vk_Nt_Ek,eq:res_vk_Nt_Ek} converge to the same limit as $N_T\to+\infty$, namely,
\begin{equation}\label{eq:asymp_limit_global}
    \lim_{N_T\to+\infty}\mathbf{v}^{(k)}_{N_t}
    = \lim_{N_T\to+\infty}\widetilde{\mathbf{v}}^{(k)}_{N_t}
    = \mathbf{v}^{(k)}_\infty
    \coloneqq \widetilde{V}_\infty^k\mathbf{u}_0,
\end{equation}
where
\begin{equation}\label{eq:V_infty}
    \widetilde{V}_\infty
    \coloneqq (I-\Psi)^{-1}(\Phi^m-\Psi)\Phi^m.
\end{equation}
\end{theorem}

\begin{proof}
We first establish the result for each modal component. For $n=1,\ldots,N$, define
$\widetilde{V}_{\infty,n}\coloneqq
(\phi_n^m-\psi_n)\phi_n^m/(1-\psi_n)$.
By \cref{lem:Ek_infty_mode}, we have
\begin{equation}
    \lim_{N_T\to+\infty}E^{(k)}_{N_T,n}
    =-\widetilde{V}_{\infty,n}^k.
\end{equation}
Since $|\phi_n|<1$, we also have $\phi_n^{N_t}\to 0$ as $N_T\to+\infty$. Thus, using \cref{eq:std_vk_Nt_Ek_mode,eq:res_vk_Nt_Ek_mode}, we obtain
\begin{align}
    \lim_{N_T\to+\infty}v^{(k)}_{N_t,n}
    &= -\lim_{N_T\to+\infty}E^{(k)}_{N_T,n}u_{0,n}
    = \widetilde{V}_{\infty,n}^k u_{0,n},\\
    \lim_{N_T\to+\infty}\widetilde{v}^{(k)}_{N_t,n}
    &= \left(-\lim_{N_T\to+\infty}E^{(1)}_{N_T,n}\right)^k u_{0,n}
    = \widetilde{V}_{\infty,n}^k u_{0,n}.
\end{align}
Hence, the MGRIT and R-MGRIT final-point approximations converge to the same limit in every mode.

To recover the full approximate solution vectors, recall that the columns $\{\mathbf{q}_n\}_{n=1}^N$ of $Q$ form an orthonormal basis. Since $N$ is finite, the modal limits can be combined to give
\begin{equation}
    \lim_{N_T\to+\infty}\mathbf{v}^{(k)}_{N_t}
    = \lim_{N_T\to+\infty}\widetilde{\mathbf{v}}^{(k)}_{N_t}
    = \sum_{n=1}^N \widetilde{V}_{\infty,n}^k u_{0,n}\mathbf{q}_n.
\end{equation}
Finally, by the simultaneous diagonalization of $\Phi$ and $\Psi$,
\begin{equation}
    Q\operatorname{diag}(\widetilde{V}_{\infty,1},\ldots,
    \widetilde{V}_{\infty,N})Q^*
    = (I-\Psi)^{-1}(\Phi^m-\Psi)\Phi^m
    = \widetilde{V}_\infty.
\end{equation}
Therefore, the common limit is $\widetilde{V}_\infty^k\mathbf{u}_0$, as claimed.
\end{proof}

By the triangle inequality, we can bound the finite-time difference between the MGRIT and R-MGRIT approximate solutions by their differences with respect to their shared asymptotic limit:
\begin{equation}\label{eq:triangle_ineq_global}
    \|\mathbf{v}^{(k)}_{N_t}-\widetilde{\mathbf{v}}^{(k)}_{N_t}\| \leq \|\mathbf{v}^{(k)}_{N_t}-\mathbf{v}^{(k)}_{\infty}\|+\|\widetilde{\mathbf{v}}^{(k)}_{N_t}-\mathbf{v}^{(k)}_{\infty}\|.
\end{equation}
In our numerical experiments from \cref{sec:num_res}, we observe that these right-hand side terms are approximately zero for a sufficiently large $N_t$, which implies the left-hand side difference is also close to zero. Namely, we will infer this in \cref{fig:mp1_grad_iter_mt} from the fact that the gradient at the C-points stalls for both MGRIT and R-MGRIT, which is consistent with $\widetilde{\mathbf{v}}^{(k)}_{N_t} \approx \mathbf{v}^{(k)}_{\infty}\approx \mathbf{v}^{(k)}_{N_t}$, with similar behavior observed in the nonlinear, nonsmooth MP2 experiments in \cref{fig:mp2-1d_grad_iter_mt,fig:mp2-2d_grad_iter_mt}.

\subsection{On a truncated R-MGRIT strategy}\label{subsec:trunc_mt}

We now motivate the idea of shortening R-MGRIT trajectories, by showing that the approximate endpoint solution obtained with one MGRIT iteration depends only weakly on the length of the trajectory.

\begin{lemma}\label{lem:truncation_decay}
Let \cref{ass:operators} hold. 
Consider an R-MGRIT trajectory with fixed
initial guess $\widetilde{\mathbf{u}}^{(k)}_0$.
Let $\widetilde{\mathbf{v}}^{(k)}_{jm}$ denote its approximation at the
$j$th C-point after one MGRIT iteration. Then the limit
$\widetilde{\mathbf{v}}^{(k)}_\infty
=\lim_{j\to\infty}\widetilde{\mathbf{v}}^{(k)}_{jm}$
exists, and
\begin{equation}\label{eq:v1_jm}
    \widetilde{\mathbf{v}}^{(k)}_{jm}
    -\widetilde{\mathbf{v}}^{(k)}_\infty
    = \Psi^{j-1}
    \left(
        \widetilde{\mathbf{v}}^{(k)}_m
        -\widetilde{\mathbf{v}}^{(k)}_\infty
    \right),
    \qquad j\geq 1.
\end{equation}
Consequently,
\begin{equation}\label{eq:asymp_trunc_err}
    \left\|
        \widetilde{\mathbf{v}}^{(k)}_{jm}
        -\widetilde{\mathbf{v}}^{(k)}_\infty
    \right\|_2
    \leq
    \rho(\Psi)^{j-1}
    \left\|
        \widetilde{\mathbf{v}}^{(k)}_m
        -\widetilde{\mathbf{v}}^{(k)}_\infty
    \right\|_2.
\end{equation}
\end{lemma}

\begin{proof}
Under constant initialization, \cref{lem:Cpt_err} gives
\begin{equation*}
    \widetilde{\mathbf{e}}^{(k)}_{jm}
    =
    \sum_{\ell=0}^{j-2}
    \Psi^{j-\ell-2}(\Phi^m-\Psi)\Phi^m
    (\Phi^{\ell m}-I)\widetilde{\mathbf{u}}^{(k)}_0,
    \qquad j\geq 2.
\end{equation*}
Separating the final term yields the recurrence
\begin{equation*}
    \widetilde{\mathbf{e}}^{(k)}_{jm}
    = \Psi\widetilde{\mathbf{e}}^{(k)}_{(j-1)m}
    +(\Phi^m-\Psi)\Phi^m
    (\Phi^{(j-2)m}-I)\widetilde{\mathbf{u}}^{(k)}_0.
\end{equation*}
Since
$\widetilde{\mathbf{v}}^{(k)}_{jm}
=\Phi^{jm}\widetilde{\mathbf{u}}^{(k)}_0
-\widetilde{\mathbf{e}}^{(k)}_{jm}$,
we obtain
\begin{align*}
    \widetilde{\mathbf{v}}^{(k)}_{jm}
    &=
    \Phi^{jm}\widetilde{\mathbf{u}}^{(k)}_0
    -\Psi\left(
        \Phi^{(j-1)m}\widetilde{\mathbf{u}}^{(k)}_0
        -\widetilde{\mathbf{v}}^{(k)}_{(j-1)m}
    \right)
    -(\Phi^m-\Psi)\Phi^m
    (\Phi^{(j-2)m}-I)\widetilde{\mathbf{u}}^{(k)}_0\\
    &=
    \Psi\widetilde{\mathbf{v}}^{(k)}_{(j-1)m}
    +(\Phi^m-\Psi)\Phi^m\widetilde{\mathbf{u}}^{(k)}_0.
\end{align*}
Using $\widetilde{\mathbf{v}}^{(k)}_m
=\Phi^m\widetilde{\mathbf{u}}^{(k)}_0$ and iterating gives
\begin{equation*}
    \widetilde{\mathbf{v}}^{(k)}_{jm}
    =
    \Psi^{j-1}\Phi^m\widetilde{\mathbf{u}}^{(k)}_0
    +\sum_{\ell=0}^{j-2}
    \Psi^\ell(\Phi^m-\Psi)\Phi^m
    \widetilde{\mathbf{u}}^{(k)}_0.
\end{equation*}
Since $\rho(\Psi)<1$, the Neumann series converges, giving
\begin{equation*}
    \widetilde{\mathbf{v}}^{(k)}_\infty
    =
    (I-\Psi)^{-1}(\Phi^m-\Psi)\Phi^m
    \widetilde{\mathbf{u}}^{(k)}_0
    =\widetilde{V}_\infty\widetilde{\mathbf{u}}^{(k)}_0.
\end{equation*}
Subtracting this limit from the preceding finite sum and using
the geometric series identity gives \cref{eq:v1_jm}.
Finally, since $\Psi$ is unitarily diagonalizable,
$\|\Psi\|_2=\rho(\Psi)<1$, which yields
\cref{eq:asymp_trunc_err}.
\end{proof}

\Cref{lem:truncation_decay} provides a motivation for shortening R-MGRIT trajectories. Although the approximate solution along a single trajectory generally converges to \(\widetilde{\mathbf{v}}^{(k)}_\infty \neq \mathbf{0}\), its difference from this limit decays geometrically at a rate governed by $\rho(\Psi)$. Thus, sufficiently far along a trajectory, computing additional C-points has little effect on its approximate endpoint. Since R-MGRIT uses only this endpoint to initialize the next trajectory, much of this additional computation may be unnecessary.

More precisely, let $\widetilde{V}_{jm}$ denote the endpoint operator for a trajectory with $j$ coarse-grid intervals. The lemma implies that $\widetilde{V}_{jm}\to\widetilde{V}_\infty$ as $j\to\infty$, and hence $(\widetilde{V}_{jm})^k\mathbf{u}_0\to\widetilde{V}_\infty^k\mathbf{u}_0=\mathbf{v}^{(k)}_\infty$ for fixed $k$. Therefore, sufficiently long trajectories can be shortened while retaining approximately the same final-point solution after $k$ restarts, potentially reducing the computational cost of R-MGRIT.

Motivated by this observation, we consider truncating each R-MGRIT trajectory to $\overline{N_t}=\gamma N_t$ fine-grid intervals, where $\gamma\in(0,1)$ and $\overline{N_t}$ is rounded to a multiple of $m$ when necessary. For sufficiently large $\overline{N_t}$, the preceding analysis suggests that, for fixed $k$,
\begin{equation}\label{eq:mgrit_asymp_limits}
    \widetilde{\mathbf{v}}^{(k)}_{\overline{N_t}}
    \approx
    \widetilde{\mathbf{v}}^{(k)}_{N_t}
    \approx
    \mathbf{v}^{(k)}_\infty.
\end{equation}
The theoretical result does not specify how short the trajectories can be made, nor does it guarantee that truncation preserves the number of restarts required for convergence. We investigate these questions numerically in \cref{sec:num_res}, using $\gamma=0.2$.

\section{Numerical results for MGRIT vs. R-MGRIT}\label{sec:num_res}

    For our numerical experiments we consider the two model problems, MP1, \cref{eq:mp1}, and MP2, \cref{eq:mp2}, with the optimization methods described in \cref{subsec:mp1,subsec:mp2} for MP1 and MP2, respectively.
    For each problem, we first determine $N_t$ via a sequential fine-grid solve, then construct our fine-grid mesh of $N_t+1$ points with spacing $\Delta t=1/L$. We choose $10^{-8}$ as the relative tolerance for all iterative methods.
    
    For R-MGRIT, as proposed in \cref{subsec:trunc_mt}, we consider truncating the fine-grid time domain by a fixed factor $\gamma=0.2$. In this case, we solve all-at-once systems of the form \cref{eq:upd_Au=g} corresponding to a fifth of the number of iterations of the underlying optimization method. 
    We illustrate convergence results for a particular set of choices for the parameters $N$ (dimension of the Laplacian matrix), $m$ (coarsening factor) and $\widehat{\ell}$ (number of levels). Later, we present results for different values of $m$, $\widehat{\ell}$ and $N$ to analyze how performance is affected by the choices of these parameters. 
    
    The initial MGRIT guess $\mathbf{u}_0 \in \mathbb{R}^{N}$ is chosen to be random with components in $[0,1]$, after which the condition $\mathbf{u}_0=\mathbf{w}_0$ is enforced on the all-at-once system \cref{eq:Au=g}.
    For both model problems, we fix a step size $s^{(\ell)}=m^{\ell-1}/L$ for the $\ell$th level operator, $\ell=1,\ldots,\widehat{\ell}$, where $m$ denotes the coarsening factor and $\widehat{\ell}$ the maximum number of levels; note that we use the same operator $\Psi$ (with the level-appropriate step size $s^{(\ell)}$) on all levels $\ell \geq 2$.

    \subsection{Numerical experiments for gradient descent applied to MP1}\label{subsec:mgrit_mp1_num_res_mt}
    
    We solve MP1, \cref{eq:mp1}, using gradient descent, \cref{eq:gd_mp1}. For this example, the 1D Laplacian matrix $A_1 \in \mathbb{R}^{N \times N}$ in \cref{eq:1d_laplacian} with 2-norm $L$ is taken with $N=n=80$, the coarsening factor is chosen to be $m=64$ and the number of levels to be $\widehat{\ell}=2$. We also choose $\mathbf{b} \in \mathbb{R}^N$ to be random with components in $[0,1]$. 
    First, we solve the problem sequentially, with $30,016$ iterations required for the convergence of the gradient descent algorithm \cref{eq:gd_mp1}. 
    %However, for $N_T=N_t/m$ to be an integer, such that the final time point is a C-point, for our MGRIT solves we work with $N_t+1=30,081$ and $N_T=470$ fine- and coarse-grid intervals, respectively. 
    MGRIT then solves \cref{eq:Au=g} using \cref{alg:mgrit}. We then consider truncating the total number of fine-grid intervals from $N_t$ to $\overline{N_t}=\gamma N_t$. For a truncation factor $\gamma=0.2$, we have $\overline{N_t}=6,080$.
    %(extending from $\gamma N_t=6,016$ to the next multiple of $m=64$) fine-grid intervals and $\overline{N_T}=95$ coarse-grid intervals. 
    %

\begin{figure}[h!]
    \centering
    \subfloat[Sequential convergence\label{fig:mp1_conv_seq_mt}]{%
        \includegraphics[width=0.45\textwidth]{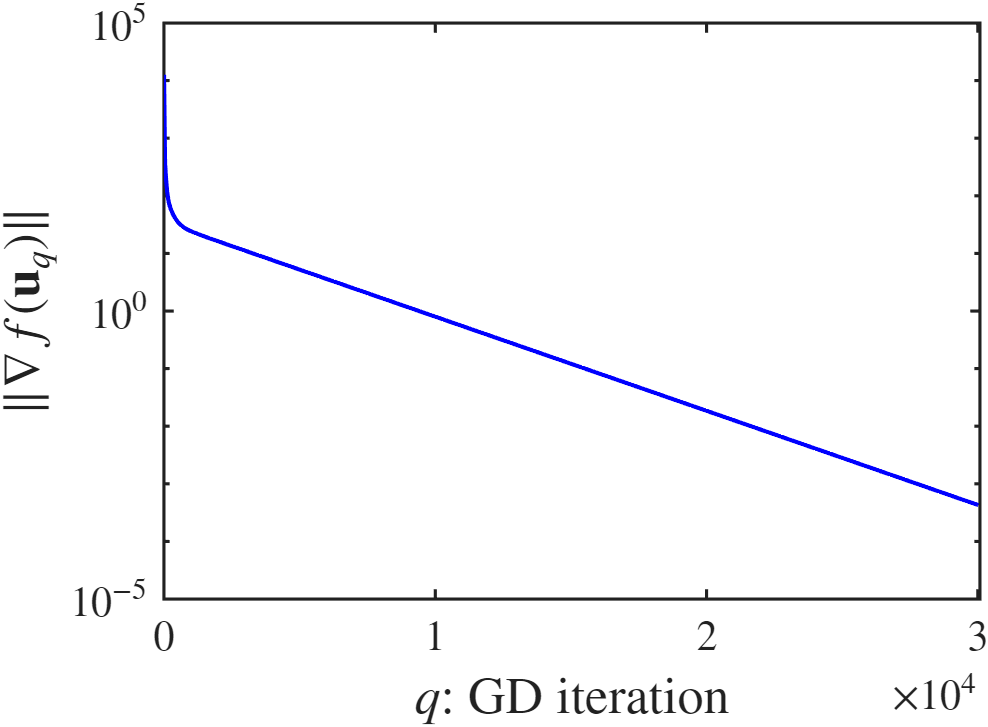}%
    }\hfill
    \subfloat[MGRIT convergence]{%
        \includegraphics[width=0.45\textwidth]{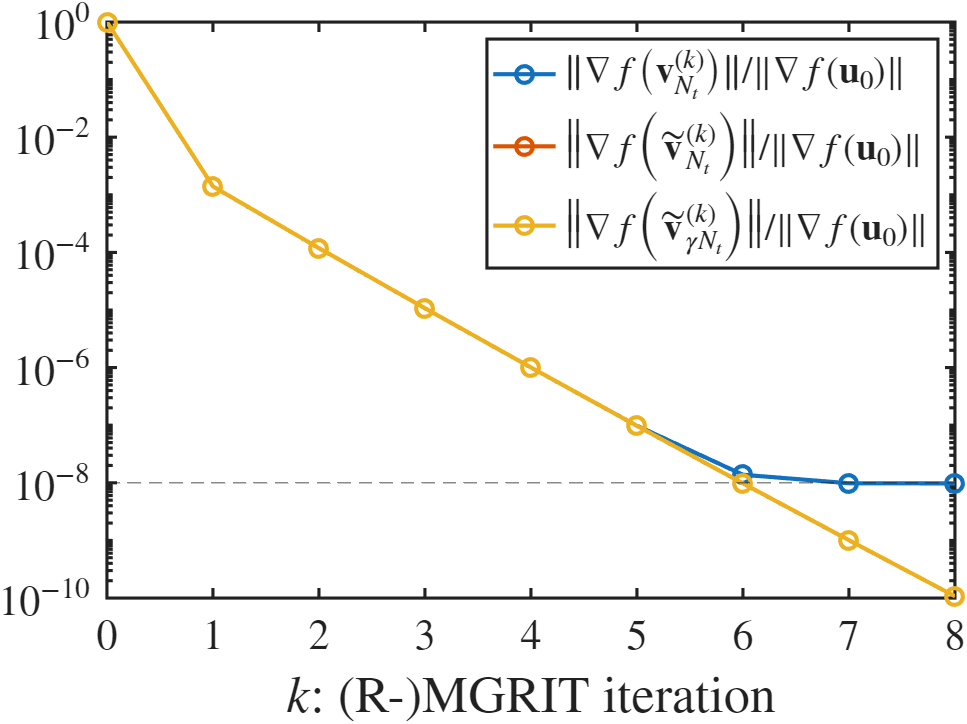}%
    }
    \caption{\textbf{(a)}: Convergence of MP1 with a standard sequential solve using gradient descent \cref{eq:gd_mp1}. 
    \textbf{(b)}: Convergence of MP1 using MGRIT (in blue), R-MGRIT (in red) and $\gamma$-truncated R-MGRIT with $\gamma=0.2$ (in yellow).}
    \label{fig:mp1_conv_mt}
\end{figure}

\begin{figure}[h!]		
    \centering
        \centering		\includegraphics[width=0.9\linewidth,trim=0.1cm 0cm 0.1cm 0cm, clip]{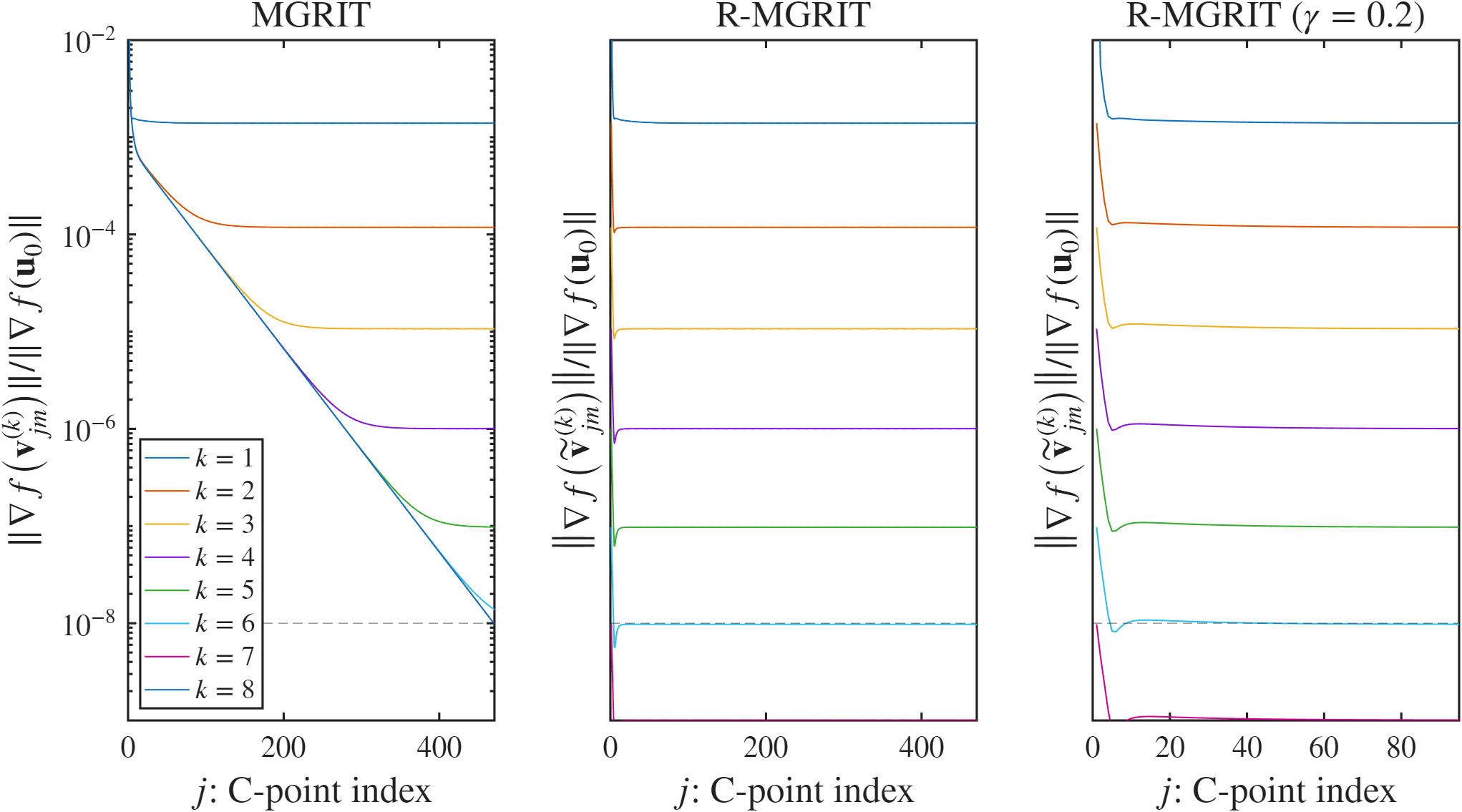}
        \vspace{-0.2cm}
         \caption{2-norm of the gradient of the MP1 approximate solutions at coarse time index $j$ and iteration $k$ of MGRIT, R-MGRIT and $\gamma$-truncated R-MGRIT (with $\gamma=0.2$), from left to right, divided by the 2-norm of the gradient of the initial guess $\mathbf{u}_0$.}	\label{fig:mp1_grad_iter_mt}
	\end{figure}

    Convergence for the sequential and MGRIT solves is illustrated in \cref{fig:mp1_conv_mt}, in which we observe that MGRIT requires 8 iterations for convergence, while R-MGRIT only requires 6 trajectories, so a total of 6 MGRIT iterations. We plot the gradient at the C-point approximate solutions of each iteration of MGRIT and R-MGRIT in \cref{fig:mp1_grad_iter_mt}, in which we observe that the truncated R-MGRIT possesses a similar convergence behavior as the full-time domain R-MGRIT solve. Therefore, by truncating $N_t$ to a sufficiently large (yet relatively small) $\overline{N_t}=\gamma N_t$, we obtain convergence in at most as many iterations as MGRIT.
    
    \Cref{fig:mp1_conv_mt}(b) shows that MGRIT and R-MGRIT converge essentially with the same speed and the cost per iteration is the same, so there is no direct benefit in considering multiple inexact trajectories instead of converging MGRIT accurately on a single trajectory. However, the middle panel of \cref{fig:mp1_grad_iter_mt} suggests that running R-MGRIT with shorter trajectories may provide a substantial cost benefit. This approach is explained in the $\gamma$-truncated numerical results in this section, and the gain in computational cost is discussed in \cref{sec:speedup}.

    Although the numerical results illustrated in \cref{fig:mp1_grad_iter_mt} are only presented for the norm of the gradient, the theory from \cref{sec:asympt_equiv} allows us to interpret the results in terms of the approximate solutions. 
    \Cref{thm:asymptotic_equivalence} predicts that the difference between the MGRIT and R-MGRIT final-point approximate solutions for each fixed iteration index $1 \leq k \leq 6$ will be small as the number of C-points becomes sufficiently large for both to approximate their common asymptotic limit of $\mathbf{v}^{(k)}_\infty$, and \cref{fig:mp1_grad_iter_mt} indicates that this happens much earlier than $N_T$. Indeed, \cref{fig:mp1_grad_iter_mt} suggests that
    \cref{eq:mgrit_asymp_limits} holds for $\overline{N_T}=\gamma N_T$ with $\gamma \leq 0.2$ with $\gamma$-truncated R-MGRIT attaining a similar gradient reduction factor to MGRIT and full-time R-MGRIT.
    
    For $k \geq 7$, we observe that $N_T$ is no longer large enough for the $\mathbf{v}^{(k)}_{N_t} \approx \mathbf{v}^{(k)}_{\infty}$ approximation to hold and the MGRIT final-point gradient stalls, with a final gradient reduction factor of approximately $10^{-8}$ as predicted by \cite[Lemma 5.1]{Araujo2026}. 
    By approximating new target trajectories rather than the initial one, R-MGRIT iterations of index $k \geq 7$ continue to reduce the gradient by approximately the same linear rate, as first illustrated in \cref{fig:mp1_conv_mt}. Thus, R-MGRIT does not require an accurate estimation of $N_t$, the number of iterations required for convergence of gradient descent \cref{eq:gd_mp1} with the prescribed $10^{-8}$ accuracy threshold.
    
    In \cref{table:1d_mp1_conv_mt}, we present MGRIT and R-MGRIT convergence results for MP1 \cref{eq:mp1} with different coarsening factors $m$, number of MGRIT levels $\widehat{\ell}$ and different sizes for the matrix $A_1$ in \cref{eq:1d_laplacian}; we note that all R-MGRIT results hold identically for a truncation factor $\gamma=0.2$ (not shown). Having established in \cite{Araujo2026} that the number of iterations does not change much as $\widehat{\ell}$ increases, here we only present results for the smallest $(\widehat{\ell}=2)$ and largest ($\widehat{\ell}=\ell$) possible number of levels given the number of time points $N_t$ and the coarsening factor $m$.
    
    \begin{table}[h!]
    \caption{MP1 \cref{eq:mp1} MGRIT convergence results for different values of $m$ and $\widehat{\ell}$. In parentheses, results for R-MGRIT. Here,  $\ell$ denotes the maximum number of levels for the given coarsening factor. For both $N=60$ and $N=80$, we have $\ell=8$ for $m=4$, $\ell=4$ for $m=16$, $\ell=3$ for $m=64$ and $\ell=2$ for $m=256$; for $N=80$, we have $m=1024$}.
    \label{table:1d_mp1_conv_mt}
        \vspace{-0.5cm}
    \centering
    \subfloat[$N=60$, $N_t=18,050$]{
        \begin{tabular}{l|l|l}
         & $\widehat{\ell}=2$ & $\widehat{\ell}=\ell$ \\ 
            \hline
            $m=4$   & 9(9) & 9(9) \\
            \hline
            $m=16$  &   8(7)    &  8(7)   \\
            \hline
            $m=64$  &   8(6)   &   8(6)  \\
            \hline
            $m=256$ &  7(5)   & --- \\
            \hline
        \end{tabular}
    }\hspace{2em}
    \subfloat[$N=80$, $N_t=30,016$]{
        \begin{tabular}{l|l|l}
         & $\widehat{\ell}=2$ & $\widehat{\ell}=\ell$ \\ 
            \hline
            $m=4$    & 9(9) & 9(9) \\
            \hline
            $m=16$  &  8(7)     &     8(7)\\
            \hline
            $m=64$  &  8(6)    &    8(6)   \\
            \hline
            $m=256$ &  8(6)   & --- \\
            \hline
            $m=1024$ &  7(5) & --- \\
        \end{tabular}
    }
    \end{table}
    
    \subsection{Numerical experiments for proximal gradient applied to MP2}\label{subsec:mgrit_mp2_num_res_mt}
    
    \subsubsection{Numerical results for MP2-1D}

    We now consider solving the EOP, MP2 \cref{eq:mp2}, using proximal gradient method \cref{eq:pg_mp2} with penalty parameter $\lambda=900$.    % 
    For MGRIT, we choose the alternating proximals operator \cref{eq:ap_mp2} to be the coarse-grid operator. We present one MP2-1D and one MP2-2D example, similar to those in \cite{Araujo2026}.
    For the MP2-1D example, we consider MP2 in \cref{eq:mp2} with $d=1$ and $N=256$, with MGRIT parameters $m=256$ and $\widehat{\ell}=2$. The obstacle $\boldsymbol{\phi}$ is obtained by discretizing $\phi(x)=\max \{0, \sin x\}$ for $x \in \Omega:=[0,3 \pi]$ on a uniform mesh of $N+2$ points. 

\begin{figure}[h!]
    \centering
    \subfloat[Sequential convergence]{%
        \includegraphics[width=0.45\textwidth]{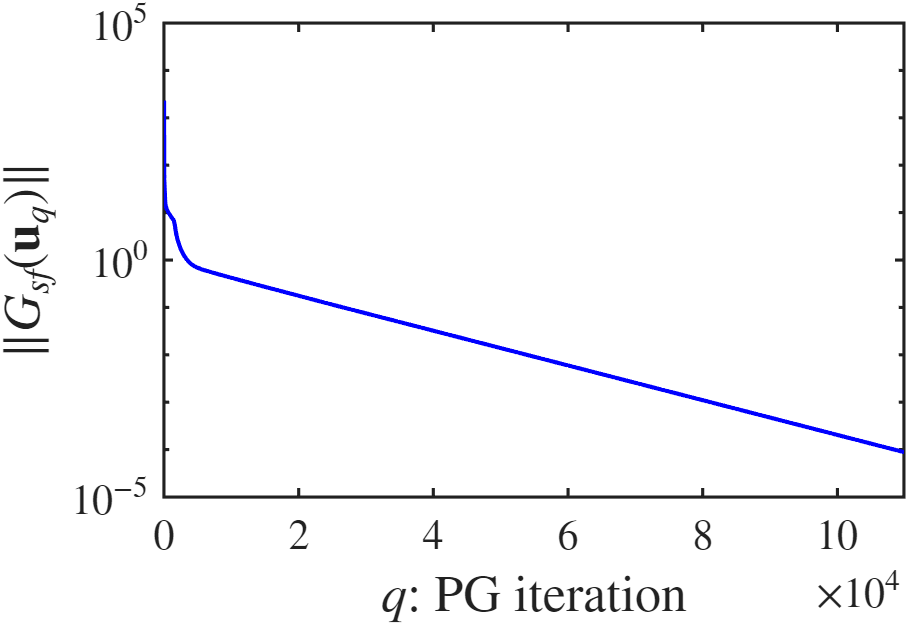}%
    }\hfill
    \subfloat[MGRIT convergence]{%
        \includegraphics[width=0.45\textwidth]{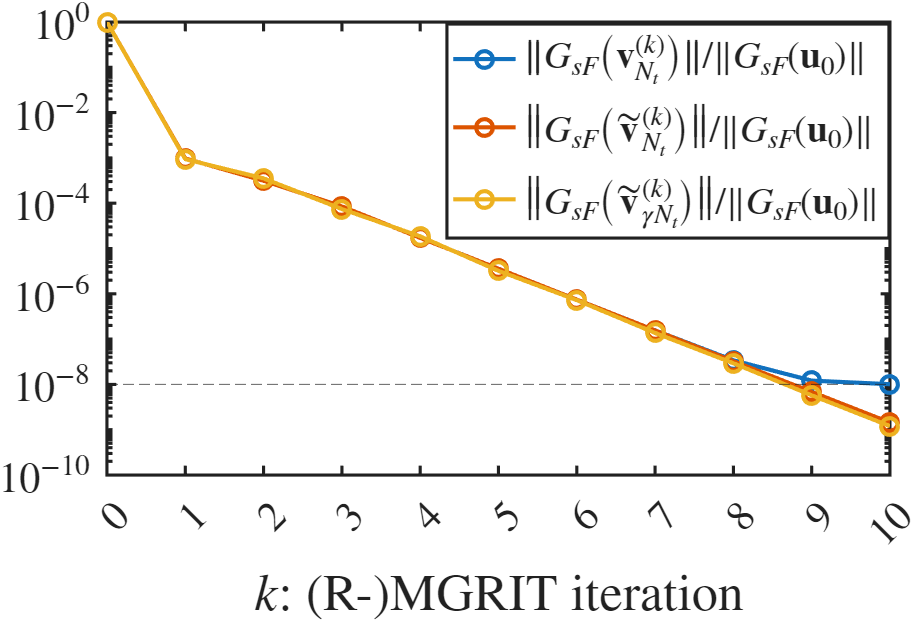}%
    }
    \caption{\textbf{(a)}: Convergence of MP2-1D with a standard sequential solve using proximal gradient \cref{eq:pg_mp2}. \textbf{(b)}: Convergence of MP2-1D using MGRIT (in blue), R-MGRIT (in red) and $\gamma$-truncated R-MGRIT with $\gamma=0.2$ (in yellow).}
    \label{fig:mp2_1d_conv_mt}
\end{figure}

\begin{figure}[h!]		
    \centering
        \centering	\includegraphics[width=0.9\linewidth,trim=0.1cm 0cm 0.1cm 0cm, clip]{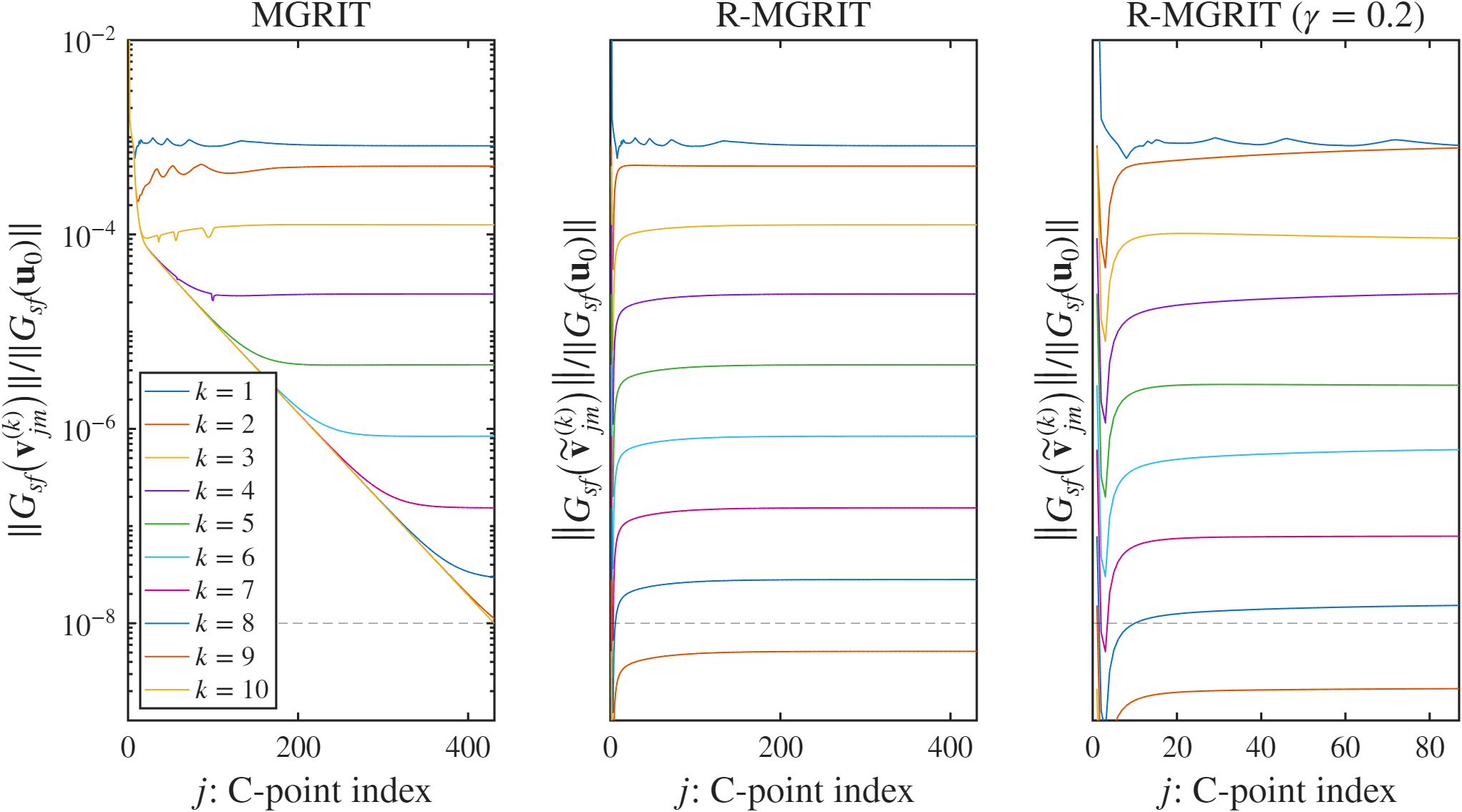}
        \vspace{-0.2cm}
         \caption{2-norm of the generalized gradient $G_{sF}$ \cref{eq:G_sF} of the MP2-1D, \cref{eq:mp2}, approximate solutions at coarse time index $j$ and iteration $k$ of MGRIT, R-MGRIT and $\gamma$-truncated R-MGRIT (with $\gamma=0.2$), from left to right, divided by the 2-norm of the gradient of the initial guess $\mathbf{u}_0$.}	\label{fig:mp2-1d_grad_iter_mt}
	\end{figure}

    We again start by presenting results for the sequential and (R-)MGRIT convergence in \cref{fig:mp2_1d_conv_mt}. 
    In this example, the proximal gradient method requires $109,888$ iterations to reach the prescribed halting tolerance, i.e., that the generalized gradient of $F$, $G_{sF}$ in \cref{eq:G_sF}, be reduced by 8 orders of magnitude; we extend the number of time points to $N_t=110,081$ such that $N_T=N_t/m=431$ is an integer. 
    In \cref{fig:mp2_1d_conv_mt}, we observe that MGRIT reaches the prescribed residual tolerance in 10 iterations, while R-MGRIT reaches it in 9 iterations both when using the full $N_t$ number of time points and when truncating the time domain to a fifth of the number of time points. The convergence speed of all 3 methods in \cref{fig:mp2_1d_conv_mt} is very similar. We present the C-point gradients of each iteration in \cref{fig:mp2-1d_grad_iter_mt}, in which we observe similar behavior in the gradient norm to what we observed in \cref{fig:mp1_grad_iter_mt} for MP1.

\subsubsection{Numerical results for MP2-2D}

\begin{figure}[h!]
    \centering
    \subfloat[Sequential convergence]{%
        \includegraphics[width=0.45\textwidth]{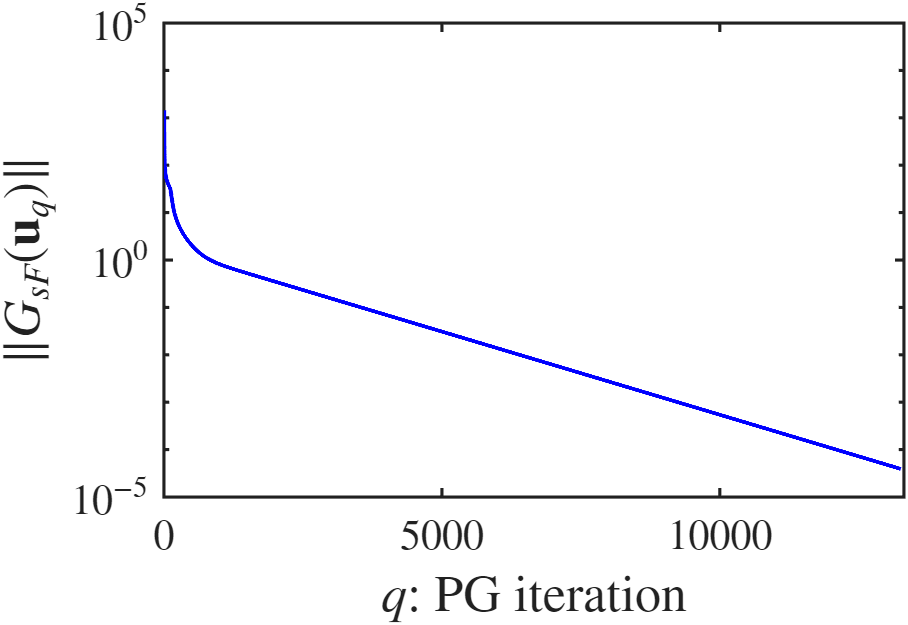}%
    }\hfill
    \subfloat[MGRIT convergence]{%
        \includegraphics[width=0.45\textwidth]{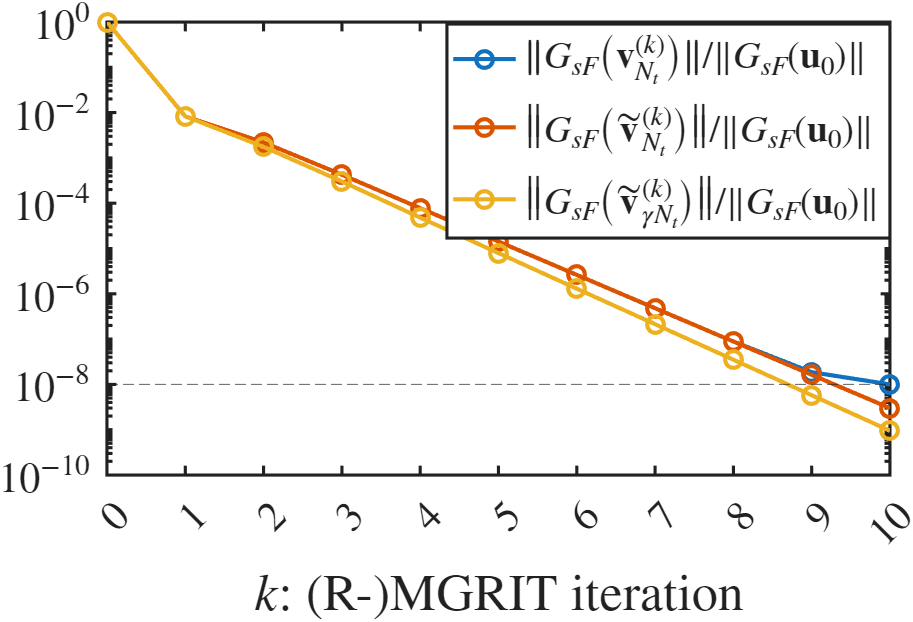}%
    }
    \caption{\textbf{(a)}: Convergence of MP2-2D with a standard sequential solve using proximal gradient \cref{eq:pg_mp2}. \textbf{(b)}: Convergence of MP2-2D using MGRIT (in blue), R-MGRIT (in red) and $\gamma$-truncated R-MGRIT with $\gamma=0.2$ (in yellow).}
    \label{fig:mp2_2d_conv_mt}
\end{figure}

\begin{figure}[h!]		
    \centering
        \centering	\includegraphics[width=0.9\linewidth,trim=0.1cm 0cm 0.1cm 0cm, clip]{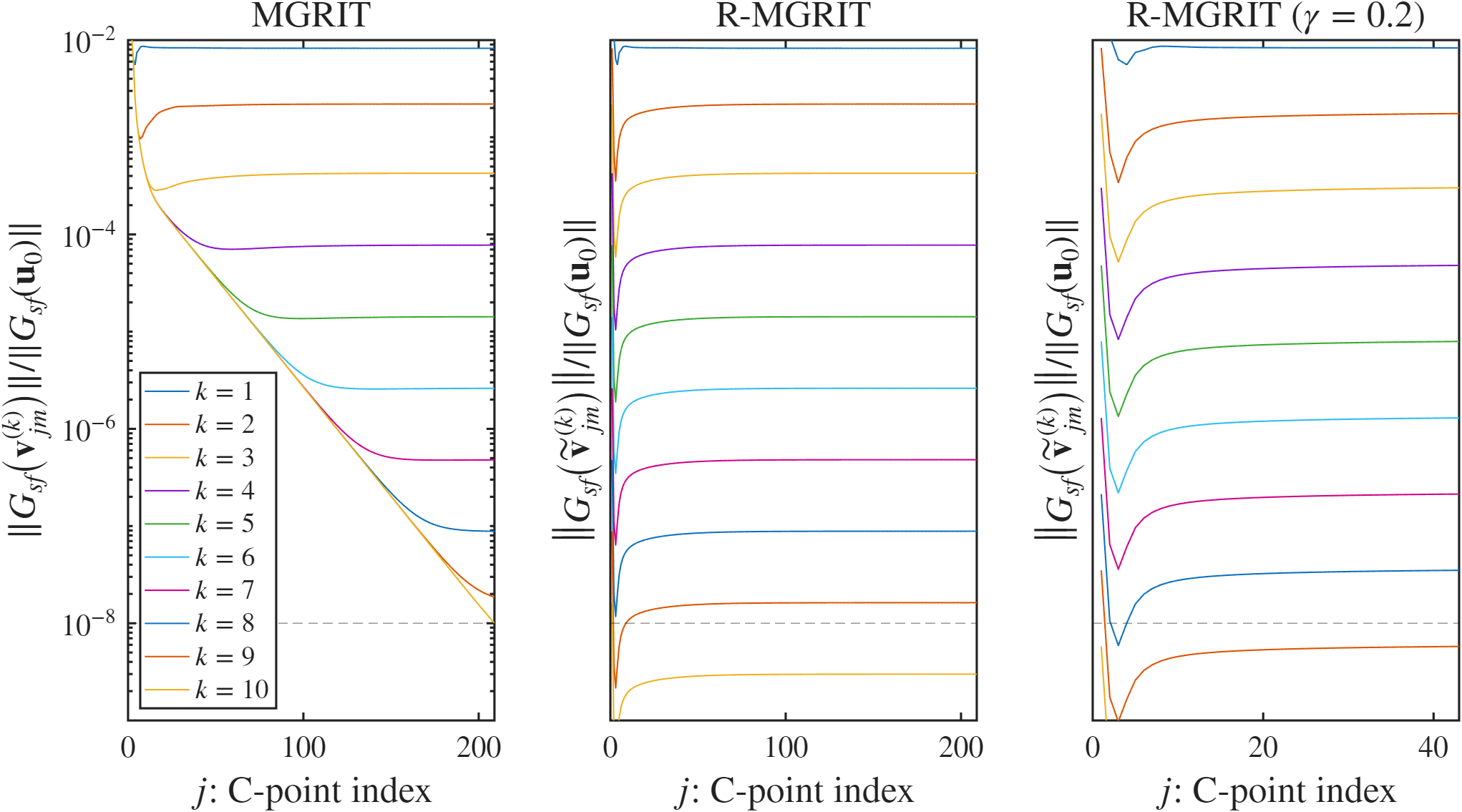}
                 \vspace{-0.2cm}
         \caption{2-norm of the generalized gradient $G_{sF}$ \cref{eq:G_sF} of the MP2-2D approximate solutions at coarse time index $j$ and iteration $k$ of MGRIT, R-MGRIT and $\gamma$-truncated R-MGRIT (with $\gamma=0.2$), from left to right, divided by the 2-norm of the gradient of the initial guess $\mathbf{u}_0$.}	\label{fig:mp2-2d_grad_iter_mt}
	\end{figure}

    For the MP2-2D example, we consider MP2 in \cref{eq:mp2} with $d=2$ and $N^2=64^2$, with MGRIT parameters $m=64$ and $\widehat{\ell}=2$. The obstacle $\boldsymbol{\phi}$ is obtained by discretizing $\phi(x, y)=\max \{0, \sin x\} \max \{0, \sin y\}$ for $(x, y) \in \Omega:=[0,3 \pi]^2$ on a uniform mesh of $N+2$ points in each spatial dimension.
    We present the sequential and MGRIT convergence results in \cref{fig:mp2_2d_conv_mt}, noting that the sequential solver required $13,256$ iterations to converge; for a coarsening factor $m=64$, we then work with $N_t=13,312$, such that $N_T=N_t/m=208$ denotes the number of coarse time intervals. We present C-point gradients in \cref{fig:mp2-2d_grad_iter_mt}, in which we observe that $\gamma$-truncated R-MGRIT (with $\gamma=0.2)$ converged in $9$ iterations, while MGRIT and full-time R-MGRIT converged in $10$ iterations. As in the MP2-1D example, we observe that $\gamma$-truncated R-MGRIT required fewer iterations to reduce the gradient by a factor of $10^{-8}$ than MGRIT and full-time R-MGRIT.

    While the theory from \cref{sec:asympt_equiv} is only presented for linear MGRIT problems, we observe generally similar results for the nonlinear MP2, \cref{eq:mp2}, in both one and two dimensions. In particular, in both \cref{fig:mp2_1d_conv_mt,fig:mp2_2d_conv_mt} we observe that R-MGRIT reduces the gradient by a factor of $10^{-8}$ in the same 9 iterations required for MGRIT to achieve a similar reduction factor, with R-MGRIT continuing to reduce the gradient in subsequent iterations while MGRIT's final-point generalized gradient stalls. Although the approximation from \cref{eq:mgrit_asymp_limits} for $\overline{N_T}=\gamma N_T$  with $\gamma=0.2$ is observed to be less accurate in that the convergence to the final gradient is slower in \cref{fig:mp2-1d_grad_iter_mt,fig:mp2-2d_grad_iter_mt} than in \cref{fig:mp1_grad_iter_mt}, as can be inferred from the rightmost panel of the respective figures, this is advantageous for $\gamma$-truncated R-MGRIT, which is able to reduce the final-point generalized gradient \cref{eq:G_sF} at a slightly faster rate because the gradient undershoots with a slow increase to the final value.

    In \cref{table:1d_mp2_conv_mt,table:2d_mp2_conv_mt}, we present convergence results for MP2, \cref{eq:mp2}, with different coarsening factors $m$, number of MGRIT levels $\widehat{\ell}$ and different sizes for the matrices $A_1$ in \cref{eq:1d_laplacian} for MP2-1D and $A_2$ in \cref{eq:2d_laplacian} for MP2-2D.
As in \cref{table:1d_mp1_conv_mt}, we only present results for the smallest $(\widehat{\ell}=2)$ and largest ($\widehat{\ell}=\ell$) permissible number of levels. 
Comparing results for MGRIT, R-MGRIT with full time domain, and $\gamma$-truncated R-MGRIT, we observe that the truncated variant achieves comparable gradient reduction in a similar number of iterations.

\begin{table}[h!]
    \caption{MP2-1D MGRIT convergence results for different values of $N$, $m$ and $\widehat{\ell}$. In parentheses, results for R-MGRIT; in square brackets, results for R-MGRIT with truncation factor $\gamma=0.2$. For $m=4$, we have $\ell=8$ for $N=256$ and $\ell=10$ for $N=512$. For both $N=256$ and $N=512$, we have $\ell=5$ for $m=16$, $\ell=4$ for $m=64$, $\ell=3$ for $m=256$ and $\ell=2$ for $m=1024$.}
    \label{table:1d_mp2_conv_mt}
        \vspace{-0.5cm}
    \centering
    \subfloat[$N=256$, $N_t=109,888$.]{
        \begin{tabular}{l|l|l}
        & $\widehat{\ell}=2$ & $\widehat{\ell}=\ell$ \\
        \hline
        $m=4$   & 10(9)[10]        & 12(13)[11] \\
        \hline
        $m=16$  & 9(8)[8] 	 &  10(9)[9]\\
        \hline
        $m=64$  & 10(9)[9]        &  10(11)[9] \\
        \hline
        $m=256$ & 10(9)[9]        & 10(9)[9]   \\
        \hline
        $m=1024$ & 11(10)[9]        &  ---        \\
        \end{tabular}
    }\hspace{2em}
    \subfloat[$N=512$, $N_t=369,612$.]{
        \begin{tabular}{l|l|l}
        & $\widehat{\ell}=2$ & $\widehat{\ell}=\ell$ \\
        \hline
        $m=4$   & 10(8)[8]   & 12(9)[10]\\
        \hline
        $m=16$  & 11(9)[9] 	 & 11(10)[9] \\
        \hline
        $m=64$  & 9(7)[9]   & 10(9)[9]  \\
        \hline
        $m=256$ & 9(8)[7]  &  --- \\
        \hline
        $m=1024$ & 10(9)[9]  &   ---       \\
        \end{tabular}
    }
\end{table}

\begin{table}[h!]
    \caption{MP2-2D MGRIT convergence results for different values of $N$, $m$ and $\widehat{\ell}$. In parentheses, results for R-MGRIT; in square brackets, results for R-MGRIT with truncation factor $\gamma=0.2$. For both $N^2=64^2$ and $N^2=128^2$, we have $\ell=7$ for $m=4$, $\ell=4$ for $m=16$, $\ell=3$ for $m=64$ and $\ell=2$ for $m=256$.}
    \label{table:2d_mp2_conv_mt}
    \vspace{-0.5cm}
    \centering
    \subfloat[$N^2=64^2$, $N_t=13,256$.]{
        \begin{tabular}{l|l|l}
        & $\widehat{\ell}=2$ & $\widehat{\ell}=\ell$ \\
        \hline
        $m=4$   & 10(10)[10] & 14(13)[13]\\
        \hline
        $m=16$  & 11(10)[10] & 12(12)[11]\\
        \hline
        $m=64$  & 10(10)[9] & 10(10)[9] \\
        \hline
        $m=256$ & 9(9)[7] & --- \\
        \end{tabular}
    }\hspace{2em}
    \subfloat[$N^2=128^2$, $N_t=46,431$.]{
        \begin{tabular}{l|l|l}
        & $\widehat{\ell}=2$ & $\widehat{\ell}=\ell$ \\
        \hline
        $m=4$   & 9(9)[9]   & 11(11)[10]\\
        \hline
        $m=16$  & 9(8)[8] 	 & 9(9)[9] \\
        \hline
        $m=64$  & 9(8)[8]   &  9(9)[8]\\
        \hline
        $m=256$ & 9(9)[8]  & ---  \\
        \end{tabular}
    }
\end{table}

\section{R-MGRIT parallel speedup estimates}\label{sec:speedup}
In \cref{sec:num_res}, we empirically observed that $\gamma$-truncated R-MGRIT achieves gradient reduction comparable to MGRIT (in the same or fewer number of iterations) with the R-MGRIT time domain truncated to a fifth of the number of time points required for MGRIT convergence to the same accuracy threshold. Because the primary motivation for employing MGRIT is to reduce the wall-clock time required to solve the sequential problem, we now formalize how this temporal truncation translates into practical parallel speedup. In the following, we assume that one R-MGRIT iteration has the same computational cost as one MGRIT iteration as described by the speedup formulas derived in \cite{Araujo2026}.
%; that is, we assume that the additional wall-clock time required for initializing each R-MGRIT all-at-once system of the form \cref{eq:upd_Au=g}, as well as any other additional operations, is negligible relative to the time required for the relaxation and coarse-grid correction steps. 
Thus, in this section we discuss 2-level R-MGRIT with a $\gamma$-truncated temporal domain, with $\gamma=1$ denoting the non-truncated R-MGRIT (or equivalently, MGRIT) speedup estimates.

\subsection{2-level $\gamma$-truncated R-MGRIT parallel speedup}\label{subsec:speedup_2_mt}

In the following, recall that $N_t$ denotes the number of time intervals on the fine-grid level $\ell=1$, and let $t_f$  denote the wall-clock time required to apply the fine-grid time-stepping operator $\Phi$ once. Similarly, recall that $N_T$ denotes the total number of time intervals on the coarse-grid level $\ell=2$, and let $t_c$ denote the wall-clock time required to apply the coarse-grid time-stepping operator $\Psi$ once.
    Let $m$ denote the coarsening factor, given by $m=N_t/N_T$. Let $N_p$ denote the number of processors. By full parallelization we mean that $N_p=N_T=N_t/m$, in which case each coarse-grid step is allocated to its own processor, along with $m$ fine-grid steps. Let $N_{it}$ denote the number of iterations required for convergence of standard MGRIT given the coarsening factor $m$. Finally, we let $\alpha \coloneqq t_c/t_f$ denote the ratio of wall-clock time required to compute one coarse-grid step relative to one fine-grid step.
    
    Given the 2-level parallel speedup estimate formula from \cite{Araujo2026}, which is based on similar formulas \cite{Minion2010, Ruprecht2012} for the parareal method \cite{Lions2001},  the speedup formula for this truncated time domain is
\begin{equation}\label{eq:s2_trunc}
    \overline{S}_2(m) = \frac{N_t}{\overline{N_{it}} \left( \alpha \overline{N_T} + 2m \right)},
\end{equation}
where $\overline{N_{it}}$ is the number of iterations required for R-MGRIT convergence on a fine grid with $\overline{N_t}=\overline{N_T}m$ intervals. Note that $\overline{N_{it}}$ still equals both the number of R-MGRIT trajectories and the total number of MGRIT iterations used by R-MGRIT.
Given some $\gamma \in (0,1]$, we substitute $\overline{N_T} =  \gamma N_t / m=\gamma N_T$ into \cref{eq:s2_trunc}, which yields
\begin{equation}\label{eq:s2_fp_gamma}
    \overline{S_2}(m, \gamma) = \frac{N_t}{\overline{N_{it}} \left( \frac{\gamma \alpha N_t}{m} + 2m \right)} = \frac{1}{\overline{N_{it}} \left( \frac{\gamma \alpha}{m} + \frac{2m}{N_t} \right)}.
\end{equation}
Note that if $\overline{N_{it}}=N_{it}$, this formula is identical to the 2-level MGRIT speedup formula from \cite{Araujo2026}, but with the operator cost ratio $\alpha$ scaled down by $\gamma$; letting $\gamma=1$ recovers the 2-level MGRIT formula exactly.
By setting the derivative of the denominator of \cref{eq:s2_fp_gamma} with respect to $m$ equal to zero for a fixed $\gamma$, we maximize $\overline{S_2}(m, \gamma)$ to find the optimal coarsening factor for the $\gamma$-truncated R-MGRIT solve:
\begin{equation}\label{eq:s2_m_gamma}
    \left.\frac{\partial}{\partial m}\left(\frac{\gamma \alpha}{m} + \frac{2m}{N_t}\right)\right|_{m=m^*_\gamma} = 0 \iff -\frac{\gamma\alpha}{(m^*_\gamma)^2} +\frac{2}{N_t}=0 \iff m^*_\gamma=\sqrt{\frac{\gamma \alpha N_t}{2}}.
\end{equation}
Let $m^*\coloneqq m^*_1$ denote the MGRIT optimal coarsening factor. Then, the relation between the truncated R-MGRIT and the MGRIT optimal coarsening factors satisfies
\begin{equation}\label{eq:m2_opt_trunc}
    \frac{m_{\gamma}^*}{m^*} 
    % = \frac{\sqrt{\gamma \alpha N_t/2}}{\sqrt{\alpha N_t/2}} 
    = \sqrt{\gamma}.
\end{equation}
By substituting $m^*$ and $m^*_\gamma$ into \cref{eq:s2_fp_gamma}, we find that the ratio between the $\gamma$-truncated R-MGRIT speedup and the standard MGRIT (or non-truncated R-MGRIT) speedup is given by
\begin{equation}
    \frac{\overline{S_2}(m_{\gamma}^*,\gamma)}{S_2(m^*)}=\frac{1}{\sqrt{\gamma}}(N_{it}/\overline{N_{it}}).
\end{equation}
In practice, we expect that  $\overline{N_{it}} \leq N_{it}$ for many small $\gamma$, as seen in \cref{table:1d_mp1_conv_mt,table:1d_mp2_conv_mt,table:2d_mp2_conv_mt} for $\gamma=0.2$. Therefore, the speedup is increased by a factor of $N_{it}/(\overline{N_{it}}\sqrt{\gamma})$ while the number of processors required for full parallelization is also reduced proportionally to $\sqrt{\gamma}$:
\begin{equation}\label{eq:processor_ratio}
 \overline{N_p^*}/N_p^*=\frac{\gamma N_t/m_\gamma^*}{N_t/m^*}=\sqrt{\gamma}.  
\end{equation}
For example, if $\gamma=1/4$ and $\overline{N_{it}}=N_{it}$, we achieve twice the speedup with half of the number of processors required for the standard MGRIT solve. 

\subsection{3-level $\gamma$-truncated R-MGRIT parallel speedup}\label{subsec:speedup_3_mt}

We now extend the previous 2-level discussion to 3-level methods based on the 3-level speedup estimates from \cite{Araujo2026}. Following reasoning similar to the 2-level case, the full-parallelization 3-level $\gamma$-truncated speedup formula is given by
    \begin{equation}\label{eq:s3_fp_gamma}
        \overline{S_3}\left(m,\gamma\right)=\frac{1}{\overline{N_{it}}\left(\frac{\gamma \alpha}{m^2}+\frac{2m}{N_t}\left(1+\frac{3 \alpha}{2}\right)\right)}.
    \end{equation}
    By setting the derivative of the denominator in \cref{eq:s3_fp_gamma} with respect to $m$ equal to zero to maximize $\overline{S_3}(m,\gamma)$, we obtain the optimal coarsening factor $m^*_\gamma$ given by 
    \begin{equation}\label{eq:s3_m_gamma}
        m^*_\gamma=\sqrt[3]{\frac{\gamma N_t}{1/\alpha+3/2}}, \quad \text{such that} \quad \frac{m_\gamma^*}{m^*} = \sqrt[3]{\gamma},
    \end{equation}
    where $m^*\coloneqq m^*_1$ denotes the 3-level MGRIT optimal coarsening factor. As in the 2-level case, truncating R-MGRIT by a factor $\gamma$ increases the speedup and reduces the processor requirement relative to MGRIT, this time proportionally to $\sqrt[3]{\gamma}$ and $\gamma^{2/3}$, respectively:
\begin{equation}\label{eq:s3_ratio}
    \frac{\overline{S_3}(m_\gamma^*, \gamma)}{S_3(m^*)}  = \frac{1}{\sqrt[3]{\gamma}} \left(N_{it}/{\overline{N_{it}}}\right), \quad \overline{N_p^*}/N_p^*=\frac{\gamma}{\sqrt[3]{\gamma}}=\gamma^{2 / 3}.
\end{equation}
For example, if $\gamma=1/8$ and $N_{it}=\overline{N_{it}}$, the predicted speedup is doubled while the processor requirement is quartered relative to the full-time 3-level MGRIT solve.

\subsection{R-MGRIT parallel speedup estimates for MP1 and MP2}\label{subsec:speedup_est_mt}

	We now present the R-MGRIT parallel speedup estimates for MP1 and MP2 in the following tables. First, we compute the optimal $m^*$ for $\widehat{\ell}=2$ and $\widehat{\ell}=3$ by \cref{eq:s2_m_gamma,eq:s3_m_gamma}, respectively, given an empirically estimated $\alpha=t_c/t_f$. We choose $m=m^*$ as the coarsening factor, rounding it to the nearest integer if needed, for the MGRIT solves. By running the MGRIT algorithm of our choice, namely \cref{alg:mgrit} for MP1 and \cref{alg:fas-mgrit} for MP2, with $m=m^*$ as the coarsening factor, we empirically obtain $N_{it}^*$ as the number of iterations required for convergence with the theoretically predicted optimal coarsening factors, and we also estimate $\alpha$. We then use $m=m^*$ and $N_{it}=N_{it}^*$ in the speedup formulas \cref{eq:s2_fp_gamma,eq:s3_fp_gamma} for 2 and 3 levels, respectively. We then do the same for truncated R-MGRIT, for which we consider a truncation factor $\gamma=0.2$, as in \cref{sec:num_res}. Note that for R-MGRIT, the $N_{it}^*$ reported in the following tables corresponds to $\overline{N_{it}}^*$ in the discussion from \cref{sec:speedup}.

    \begin{table}[h!]
    \caption{Speedup estimates for MP1, \cref{eq:mp1}, for different values of $N$ and $\widehat{\ell}$. Here, $m^*$ and $S^*$ are computed based on the numerically estimated $\frac{t_c}{t_f}$ values given by $\alpha=3.1$ for $N=60$, and $\alpha=3.5$ for $N=80$. $N_p$ denotes the number of processors required to achieve optimal speedup.}
    \label{table:mp1_speedup_mt}
    \centering
    \begin{tabular}{p{0.15\textwidth} c c}
        MGRIT & 
        \begin{tabular}{l|l|l|l|l}
            $\widehat{\ell}$ & $m^*$ & $N_{it}^*$ & $S^*_{\widehat{\ell}}$ & $N_p$\\
            \hline
            $2$ & $167$ & $8$  & $3.37$ & $109$\\
            \hline
            $3$ & $21$ & $8$ & $6.20$ & $860$
        \end{tabular} &
        \begin{tabular}{l|l|l|l|l}
            $\widehat{\ell}$ & $m^*$ & $N_{it}^*$ & $S^*_{\widehat{\ell}}$ & $N_p$ \\
            \hline
            $2$ &  $229$ & $8$ & $4.09$ & $132$ \\
            \hline
            $3$ & $26$ & $8$ & $7.81$ & $1155$
        \end{tabular} \\[6ex]
        
        R-MGRIT \newline ($\gamma =0.2$) &
        \begin{tabular}{l|l|l|l|l}
            $\widehat{\ell}$ & $m^*$ & $N_{it}^*$ & $S^*_{\widehat{\ell}}$ & $N_p$\\
            \hline
            $2$ & $74$ & $6$  & $10.05$ & $49$ \\
            \hline
            $3$ & $13$ & $6$ & $14.12$ & $278$
        \end{tabular} &
        \begin{tabular}{l|l|l|l|l}
            $\widehat{\ell}$ & $m^*$ & $N_{it}^*$ & $S^*_{\widehat{\ell}}$ & $N_p$ \\
            \hline
            $2$ & $103$  & $6$ & $12.20$ & $59$ \\
            \hline
            $3$ & $15$ & $8$ & $13.36$ &  $401$
        \end{tabular} \\[4ex]
        & \multicolumn{1}{c}{(a) $N=60$, $N_t=18,050$} & \multicolumn{1}{c}{(b) $N=80$, $N_t=30,016$}
    \end{tabular}
\end{table}

\begin{table}[h!]
    \caption{Speedup estimates for MP2-1D, \cref{eq:mp2}, for different values of $N$ and $\widehat{\ell}$. Here, $m^*$ and $S^*$ are computed based on the numerically estimated $\frac{t_c}{t_f}$ values given by $\alpha=2.4$ for $N=256$, and $\alpha=3.9$ for $N=512$. $N_p$ denotes the number of processors required to achieve optimal speedup.}
    \label{table:1d_mp2_speedup_mt}
    \centering
    \begin{tabular}{p{0.15\textwidth} c c}
        MGRIT &
        \begin{tabular}{l|l|l|l|l}
            $\widehat{\ell}$ & $m^*$ & $N_{it}^*$ & $S^*_{\widehat{\ell}}$ & $N_p$ \\
            \hline
            $2$ & $363$ & $10$ & $7.56$ & $303$\\
            \hline
            $3$ & $39$ & $10$ & $20.65$ & $2818$
        \end{tabular} &
        \begin{tabular}{l|l|l|l|l}
            $\widehat{\ell}$ & $m^*$ & $N_{it}^*$ & $S^*_{\widehat{\ell}}$ & $N_p$  \\
            \hline
            $2$ & $849$ & $10$ & $10.88$ & $436$\\
            \hline
            $3$ & $59$ & $9$ & $33.60$ & $6265$
        \end{tabular} \\[6ex]
        
        R-MGRIT \newline ($\gamma =0.2$) &
        \begin{tabular}{l|l|l|l|l}
            $\widehat{\ell}$ & $m^*$ & $N_{it}^*$ & $S^*_{\widehat{\ell}}$ & $N_p$\\
            \hline
            $2$ & $162$ & $9$  & $18.80$ & $136$ \\
            \hline
            $3$ & $23$  & $10$  & $35.30$  & $956$
        \end{tabular} &
        \begin{tabular}{l|l|l|l|l}
            $\widehat{\ell}$ & $m^*$ & $N_{it}^*$ & $S^*_{\widehat{\ell}}$ & $N_p$ \\
            \hline
            $2$ & $380$  & $7$ & $34.79$ & $195$ \\
            \hline
            $3$ & $35$ & $8$ & $64.63$ & $2113$
        \end{tabular} \\[4ex]
        & \multicolumn{1}{c}{(a) $N=256$, $N_t=109,888$} & \multicolumn{1}{c}{(b) $N=512$, $N_t=369,612$}
    \end{tabular}
\end{table}

\begin{table}[h!]
    \caption{Speedup estimates for MP2-2D, \cref{eq:mp2}, for different values of $N$ and $\widehat{\ell}$. Here, $m^*$ and $S^*$ are computed based on the numerically estimated $\frac{t_c}{t_f}$ values given by $\alpha=13$ for $N^2=64^2$, and $\alpha=7$ for $N^2=128^2$. $N_p$ denotes the number of processors required to achieve optimal speedup.}
    \label{table:2d_mp2_speedup_mt}
    \centering
    \begin{tabular}{p{0.15\textwidth} c c}
        MGRIT &
        \begin{tabular}{l|l|l|l|l}
            $\widehat{\ell}$ & $m^*$ & $N_{it}^*$ & $S^*_{\widehat{\ell}}$ & $N_p$ \\
            \hline
            $2$ & $294$  & $9$ &  $1.25$ & $46$ \\
            \hline
            $3$ & $20$ & $10$ & $1.06$  & $663$
        \end{tabular} &
        \begin{tabular}{l|l|l|l|l}
            $\widehat{\ell}$ & $m^*$ & $N_{it}^*$ & $S^*_{\widehat{\ell}}$ & $N_p$ \\
            \hline
            $2$ & $403$  & $9$ &  $3.19$ & $116$\\
            \hline
            $3$ & $30$ & $9$  & $4.91$ & $1548$
        \end{tabular} \\[6ex]
        
        R-MGRIT \newline ($\gamma =0.2$) &
        \begin{tabular}{l|l|l|l|l}
            $\widehat{\ell}$ & $m^*$ & $N_{it}^*$ & $S^*_{\widehat{\ell}}$ & $N_p$\\
            \hline
            $2$ & $131$ &  $8$ & $3.16$ & $21$\\
            \hline
            $3$ & $12$ & $11$ & $1.65$ & $221$
        \end{tabular} &
        \begin{tabular}{l|l|l|l|l}
            $\widehat{\ell}$ & $m^*$ & $N_{it}^*$ & $S^*_{\widehat{\ell}}$ & $N_p$ \\
            \hline
            $2$ & $180$  & $8$ & $8.05$ & $52$\\
            \hline
            $3$ & $18$ & $9$ & $8.39$  & $516$
        \end{tabular} \\[4ex]
        & \multicolumn{1}{c}{(a) $N^2=64^2$, $N_t=13,256$} & \multicolumn{1}{c}{(b) $N^2=128^2$, $N_t=46,431$}
    \end{tabular}
\end{table}

In \cref{table:mp1_speedup_mt,table:1d_mp2_speedup_mt,table:2d_mp2_speedup_mt}, we observe that truncating the time domain to a fifth of the number of time points required for standard MGRIT convergence provides a substantial improvement of the standard MGRIT speedup, in accordance with the theory developed in \cref{sec:speedup}. Notably, our results here suggest that R-MGRIT converges in less parallel wall-clock time by employing only a small fraction of $N_t$ on each trajectory, the number of iterations required for the underlying optimization method to converge to a sufficiently accurate approximation of the true optimum.

We note that truncating more aggressively, e.g, choosing $\gamma=0.125$ instead of $\gamma=0.2$, did not significantly change the predicted speedup. We remark that an optimal truncation factor $\gamma$ can not be determined \textit{a priori} since we cannot estimate $N^*_{it}$ as a function of $\gamma$ and the other parameters. However, our numerical results suggest $\gamma=0.2$ is a good baseline parameter choice for our model problems.

\section{Applying Anderson-acceleration to R-MGRIT}\label{sec:anderson}

In \cref{sec:asympt_equiv,sec:num_res}, we observed that R-MGRIT is not expected to achieve faster gradient reduction than MGRIT. Thus, in an attempt to obtain faster gradient reduction over MGRIT, we consider applying acceleration techniques to R-MGRIT. Namely, we consider applying Anderson acceleration \cite{Anderson1965,DeSterck2022} to accelerate R-MGRIT.
To this end, we define $\mathcal{F}$ as the map that takes an
initial guess and returns the approximate final point obtained
after one MGRIT iteration on the corresponding optimization
trajectory. Since R-MGRIT uses this final point to initialize
the next trajectory, we can write R-MGRIT as the fixed-point
iteration
\begin{equation}\label{eq:rmgrit_fix_pt}
    \widetilde{\mathbf{u}}^{(k+1)}_{0}
    =\mathcal{F}(\widetilde{\mathbf{u}}^{(k)}_{0}),
    \quad k=0,1,\ldots.
\end{equation}
The minimizer $\mathbf{u}^*$ of the underlying optimization
problem, $\min_{\mathbf{u}\in\mathbb{R}^N}F(\mathbf{u})$,
is a fixed point of this map, i.e.,
$\mathcal{F}(\mathbf{u}^*)=\mathbf{u}^*$.
%
% Here, $\mathcal{F}:\mathbb{R}^N \to \mathbb{R}^N$ denotes the mapping applied to the initial guess of the $k$th all-at-once system, $\widetilde{\mathbf{u}}^{(k)}_{0}$, which yields the R-MGRIT final-time approximate solution of the $k$th all-at-once system, which we take to be the initial guess for $(k+1)$th all-at-once system, $\widetilde{\mathbf{u}}^{(k+1)}_{0}\coloneqq \widetilde{\mathbf{v}}^{(k)}_{N_t}$.  
% %
% In particular, for linear problems with one MGRIT iteration per restart, by \cref{eq:res_vk_Nt_Ek} we have that $\mathcal{F}(\widetilde{\mathbf{u}}^{(k)}_{0})=\widetilde{V}_{N_t}\widetilde{\mathbf{u}}^{(k-1)}_{0}$, where $\widetilde{V}_{N_t}$ is defined in \cref{eq:restarted_single_iter_op}. 

The Anderson acceleration iteration for \cref{eq:rmgrit_fix_pt} given a window size $w \in \mathbb{N} \cup \{+\infty\}$, denoted by AA($w$), is given by
\begin{equation}\label{eq:aa}
\widetilde{\mathbf{u}}^{(k+1)}_{0}=\mathcal{F}\left(\widetilde{\mathbf{u}}^{(k)}_{0}\right)+\sum_{i=1}^{\min (k, w)} \beta_i^{(k)}\left(\mathcal{F}\left(\widetilde{\mathbf{u}}^{(k)}_{0}\right)-\mathcal{F}\left(\widetilde{\mathbf{u}}^{(k-i)}_{0}\right)\right), \quad k=0,1, \ldots,
\end{equation}
where for each $k$ the coefficients $\beta_i^{(k)}$ are determined by solving 
\begin{equation}\label{eq:aa_beta}
\min _{\left\{\beta_i^{(k)}\right\}}\left\|r\left(\widetilde{\mathbf{u}}^{(k)}_{0}\right)+\sum_{i=1}^{\min (k, w)} \beta_i^{(k)}\left(r\left(\widetilde{\mathbf{u}}^{(k)}_{0}\right)-r\left(\widetilde{\mathbf{u}}^{(k-i)}_{0}\right)\right)\right\|,
\end{equation}
where $r(\widetilde{\mathbf{u}}^{(k)}_{0})\coloneqq \widetilde{\mathbf{u}}^{(k)}_{0}-\mathcal{F}(\widetilde{\mathbf{u}}^{(k)}_{0})$ for all $k$. Put simply, AA($w$) attempts to accelerate \cref{eq:rmgrit_fix_pt} by using the $w$ previous initial guesses to generate an improved initial guess for the next all-at-once system. 
In the tests below we use AA(1), in which case only the previous initial guess is used, and AA($\infty$), corresponding to using the full available history of initial guesses.
%Note that if the original fixed-point iteration \cref{eq:rmgrit_fix_pt} converges to the desired accuracy in $q$ iterations, by choosing $w \geq q$ we ensure that all previous initial guesses are used in this process.

\begin{figure}[h!]
    \centering
    \subfloat[AA($1$) applied to R-MGRIT for MP1.]{%
        \includegraphics[width=0.45\textwidth]{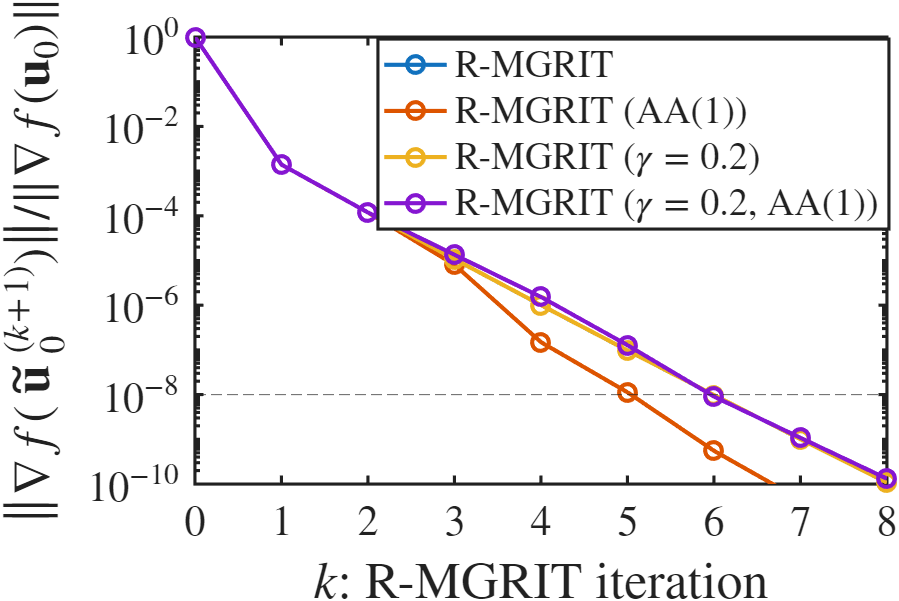}%
    }\hfill
    \subfloat[AA($1$) applied to R-MGRIT for MP2-2D.]{%
        \includegraphics[width=0.45\textwidth]{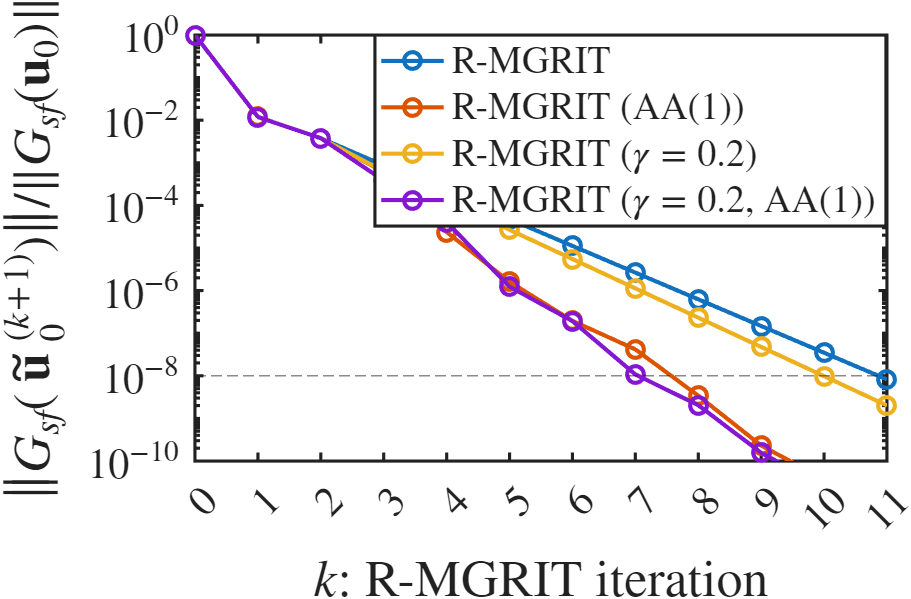}%
    }
    \\    
    \subfloat[AA($\infty$) applied to R-MGRIT for MP1.]{%
        \includegraphics[width=0.45\textwidth]{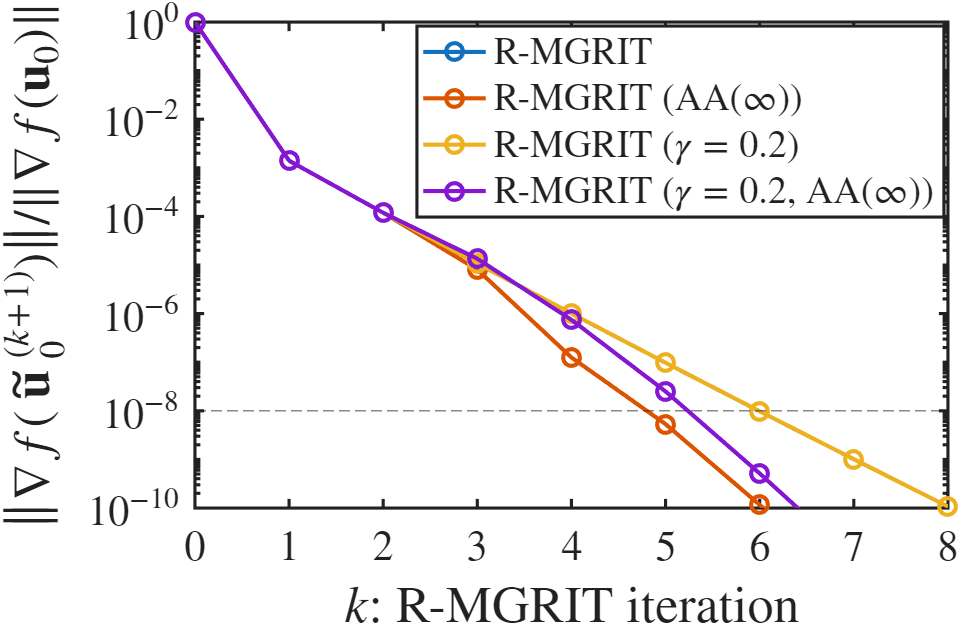}%
    }\hfill
    \subfloat[AA($\infty$) applied to R-MGRIT for MP2-2D.]{%
        \includegraphics[width=0.45\textwidth]{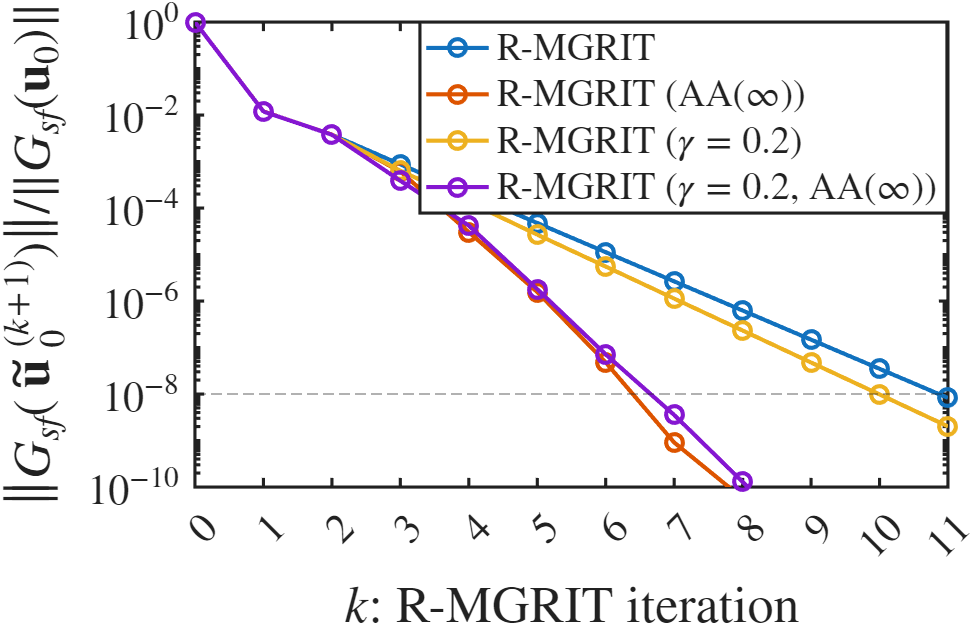}%
    }
    \caption{Convergence of R-MGRIT (blue), full-time R-MGRIT with Anderson acceleration \cref{eq:aa} (orange), $\gamma$-truncated R-MGRIT (yellow) and $\gamma$-truncated R-MGRIT with Anderson acceleration (purple) variants, with $\gamma=0.2$, for: \textbf{(a)} MP1 and AA(1), \textbf{(b)} MP2-2D and AA(1), \textbf{(c)} MP1 and AA($\infty$), and \textbf{(d)} MP2-2D and AA($\infty$). In \textbf{(a)} and \textbf{(c)}, the blue curve sits underneath the yellow curve.}
    \label{fig:aa_conv_mt}
\end{figure}

In \cref{fig:aa_conv_mt}, we present numerical experiments similar to those in \cref{sec:num_res} in which we apply Anderson acceleration to R-MGRIT with a full time domain and truncated R-MGRIT with $\gamma=0.2$.
We consider the setting of the numerical experiments illustrated in \cref{fig:mp1_conv_mt,fig:mp2_2d_conv_mt} for MP1-1D, \cref{eq:mp1}, and MP2-2D, \cref{eq:mp2}, respectively; this time, we compare MGRIT with Anderson-accelerated R-MGRIT and Anderson-accelerated truncated R-MGRIT.
We observe in \cref{fig:aa_conv_mt} that AA($\infty$) reduces the number of R-MGRIT iterations required for convergence in experiments for both model problems. AA(1) achieves a similar improvement in gradient reduction in all instances, except for the $\gamma$-truncated MP1 experiment. Since Anderson acceleration adds trivial overhead cost to the underling fixed-point iteration, the acceleration provided would improve the parallel speedup estimates presented in \cref{sec:speedup} by further reducing $\overline{N_{it}}$, the number of R-MGRIT iterations required for convergence, relative to $N_{it}$, the number of MGRIT iterations required for convergence.

\section{Conclusions}\label{sec:conclusions}

We have extended the parallel-in-iteration MGRIT framework for optimization from \cite{Araujo2026} by introducing restarted MGRIT (R-MGRIT). Rather than repeatedly improving the approximation of a single optimization trajectory, \linebreak R-MGRIT performs one MGRIT iteration per trajectory, using its approximate final point to initialize the next. For linear problems under suitable assumptions, we proved that, for any fixed $k$, the final-point approximations obtained from $k$ MGRIT iterations and $k$ R-MGRIT trajectories converge to the same limit as the trajectory length tends to infinity.

We further showed that the approximate endpoint of a single R-MGRIT trajectory approaches its infinite-trajectory limit geometrically, which together with our numerical observations motivates a truncated R-MGRIT strategy in which each trajectory is shortened by a factor $\gamma\in(0,1)$. Our numerical experiments on linear and nonlinear optimization problems showed that truncation generally does not reduce gradient reduction per trajectory, while leading to substantial performance gains: our parallel cost analysis predicts that shortening the trajectories can substantially increase parallel speedup while reducing the number of processors required.

Finally, we investigated the use of Anderson acceleration with R-MGRIT and observed further reductions in the number of trajectories required for convergence in our numerical experiments. Future work will consider extensions of these parallel-in-iteration methods to broader classes of optimization problems and optimization algorithms, along the lines discussed in \cite{Araujo2026}.

\section*{Acknowledgements}
The authors used Google's Gemini and OpenAI's ChatGPT during the final stages of the preparation of this manuscript.
All AI suggestions mathematical claims, proofs, computations, citations, and final wording were reviewed, verified, and edited by the authors, who take full responsibility for the manuscript.
Los Alamos National Laboratory report number LA-UR-26-28803.

\bibliography{references_doi}

\end{document}